\documentclass[12pt]{article}
\usepackage{amsfonts}

\usepackage{bm}
\usepackage{enumerate} 
\usepackage{amssymb,amsmath}
\usepackage{caption}
\usepackage{subcaption}
\usepackage{accents} 

\usepackage{dsfont}

\usepackage{stackengine}

\usepackage{accents}

\usepackage{tikz}
\usetikzlibrary{automata,topaths}
\usetikzlibrary{shapes}
\usetikzlibrary{plotmarks}

\usepackage{float} 

\usepackage{color} 
\definecolor{lightblue}{rgb}{0,0.2,0.5}

\usepackage[colorlinks=true, urlcolor=lightblue,linkcolor=lightblue, citecolor=lightblue]{hyperref}

\usepackage{bbm}

\DeclareMathOperator{\imag}{Im}

\DeclareMathAlphabet{\eufrak}{U}{}{}{} 
\SetMathAlphabet\eufrak{normal}{U}{euf}{m}{n}
\SetMathAlphabet\eufrak{bold}{U}{euf}{b}{n}

\newcommand{\bigabs}[1]{\left| #1 \right|}
\newcommand{\eps}{\varepsilon}

\newcommand{\R}{\mathbb{R}}

\newcommand{\T}{\mathbb{T}}
\newcommand{\C}{\mathbb{C}}
\newcommand{\E}{\mathbb{E}}

\newcommand{\biganglebracket}[1]{\left\langle #1 \right\rangle}

\newcommand{\bigbracket}[1]{\left( #1 \right)}
\newcommand{\bignorm}[1]{\left\| #1 \right\|}

\newcommand{\N}{\mathbb{N}}

\newcommand{\bigsquarebracket}[1]{\left[ #1 \right]}

\newtheorem{prop}{Proposition}[section]

\newtheorem{condition}[prop]{Condition}
\newtheorem{lemma}[prop]{Lemma}
 
\newtheorem{definition}[prop]{Definition}
\newtheorem{corollary}[prop]{Corollary}
\newtheorem{theorem}[prop]{Theorem}

\newtheorem{example}[prop]{Example}

\makeatletter
\newcommand*\rel@kern[1]{\kern#1\dimexpr\macc@kerna}
\newcommand*\widebar[1]{
	\begingroup
	\def\mathaccent##1##2{
		\rel@kern{0.8}
		\overline{\rel@kern{-0.8}\macc@nucleus\rel@kern{0.2}}
		\rel@kern{-0.2}
	}
	\macc@depth\@ne
	\let\math@bgroup\@empty \let\math@egroup\macc@set@skewchar
	\mathsurround\z@ \frozen@everymath{\mathgroup\macc@group\relax}
	\macc@set@skewchar\relax
	\let\mathaccentV\macc@nested@a
	\macc@nested@a\relax111{#1}
	\endgroup
}
\makeatother

\makeatletter
\DeclareRobustCommand\widecheck[1]{{\mathpalette\@widecheck{#1}}}
\def\@widecheck#1#2{
	\setbox\z@\hbox{\m@th$#1#2$}
	\setbox\tw@\hbox{\m@th$#1
		\widehat{
			\vrule\@width\z@\@height\ht\z@
			\vrule\@height\z@\@width\wd\z@}$}
	\dp\tw@-\ht\z@
	\@tempdima\ht\z@ \advance\@tempdima2\ht\tw@ \divide\@tempdima\thr@@
	\setbox\tw@\hbox{
		\raise\@tempdima\hbox{\scalebox{1}[-1]{\lower\@tempdima\box
				\tw@}}}
	{\ooalign{\box\tw@ \cr \box\z@}}}
\makeatother

\def\({\left(}
\def\){\right)}

\def\[{\left[}
\def\]{\right]}
\def\real{{\mathord{\mathbb R}}}
\def\N{{\mathord{\mathbb N}}}

\def\Var{\mathrm{Var}}

\def\P{\mathbb{P}}

\newcommand{\Z}{\mathbb{Z}}

\DeclareMathOperator{\re}{Re}

\newenvironment{Proof}{\removelastskip\par\medskip
	\noindent{\em Proof.} \rm}{\penalty-20\null\hfill$\square$\par\medbreak}

\newenvironment{Proofy}{\removelastskip\par\medskip
	\noindent{\em Proof} \rm}{\penalty-20\null\hfill$\square$\par\medbreak}

\allowdisplaybreaks

\numberwithin{equation}{section}

\usepackage{graphicx}
\usepackage{flushend,cuted}
\usepackage{bm}
\usepackage{tabularx}

\usepackage{indentfirst}
\usepackage{amssymb}
\usepackage{xparse}
\usepackage{tikz}
\usepackage{mdwlist}
\usepackage{tkz-graph}
\usepackage{mathdots} 

\GraphInit[vstyle = Shade]
\usetikzlibrary[intersections,
positioning,
petri,
backgrounds,
fit,
decorations.pathmorphing,
arrows,
arrows.meta,
bending,
calc,
intersections,
through,
backgrounds,
shapes.geometric,
quotes,
matrix,
trees,
shapes.symbols,
graphs,
math,
patterns,
external,
scopes,
matrix,
lindenmayersystems,
shapes.callouts,
shapes.misc,
angles,
shapes.arrows,
shadings]

\usetikzlibrary{arrows,automata,shapes,positioning,decorations.pathmorphing,snakes}

\tikzset{snake it/.style={-stealth,
		decoration={snake, 
			amplitude = .4mm,
			segment length = 2mm,
			post length=0.9mm},decorate}}

\usetikzlibrary{matrix,calc}

\newcommand{\dx}{\, dx}

\newcommand{\du}{\, du}

\renewcommand{\P}{\mathbb{P}}

\newcommand{\disp}{\displaystyle}

\usepackage{mathdots}

\usepackage{empheq} 

\usepackage{upgreek}

\tikzset{hide labels/.style={every label/.append style={text opacity=0}}}

\begin{document}
	\title{
		\huge
				Mixing rates for linear operators under infinitely divisible
		measures on Banach spaces 
	} 
	
	\author{
		Camille Mau\footnote{\href{mailto:camille001@e.ntu.edu.sg}{camille001@e.ntu.edu.sg}}
		\qquad
		Nicolas Privault\footnote{
			\href{mailto:nprivault@ntu.edu.sg}{nprivault@ntu.edu.sg}
		}
		\\
		\small
		Division of Mathematical Sciences
		\\
		\small
		School of Physical and Mathematical Sciences
		\\
		\small
		Nanyang Technological University
		\\
		\small
		21 Nanyang Link, Singapore 637371
	}
	
	\maketitle
	
	\vspace{-0.5cm}
	
	\begin{abstract} 
		We derive rates of convergence for
		the mixing of operators 
                under infinitely divisible measures
		in the framework of linear dynamics
		on Banach spaces. 
		Our approach is based on the characterization
		of mixing in terms of codifference functionals
		and control measures, and extends previous results
		obtained in the Gaussian setting via the use of
		covariance operators. 
		Explicit mixing rates are obtained for weighted
		shifts 
		under compound Poisson, $\alpha$-stable,
		and tempered $\alpha$-stable 
		measures. 
			\end{abstract}
		
	
	\baselineskip0.7cm

	\section{Introduction}
	\noindent
	The mixing and ergodicity properties of Gaussian processes
	and dynamical systems under Gaussian measures 
	have been originally studied in \cite{akutowicz}
	and \cite{cornfeld},
	see Chapter~14-\S2 and Theorems~1 and~2
	therein, 
	in connection with the spectral properties
	of unitary transformations, spectral measures
	and Gaussian covariances. 
	
	\medskip 
	
	On the other hand, characterizations of mixing 
	of 
	continuous linear operators $T:E\to E$ invariant on a complex 
	separable Banach space~$E$ have been obtained in the framework of {linear dynamics}
	under a Gaussian measure  $\mu$ on a complex Banach space~$E$. 
	Recall (see, e.g., \cite[Definition 5.23]{DLO}) that
	a measure-preserving map $T$ 
	on $(E, \mu)$ 
	is strongly mixing
	if either of the two following equivalent conditions
	is satisfied: 
	\begin{enumerate}[(i)]
		\item $\displaystyle
		\lim_{n\to \infty} \mu (A \cap  T^{-n} (B)) = \mu (A)\mu (B)$,
		$A,B\in {\cal B}$,
		\item
		$\displaystyle
		\lim_{n\to \infty} 
		I_n(f,g) = 0$, 
		$f,g \in L^2( E , \mu)$, 
	\end{enumerate}
	where ${\cal B}$ is the Borel $\sigma$-algebra of $E$ and 
	$$
	I_n(f,g) : = \int_E f(z) g(T^nz) \mu (dz)
	-      \int_E f(z)  \mu (dz)
	\int_E g(z) \mu (dz), \quad n\geq 0. 
	$$
	Likewise, $T$ is 
	weakly mixing with respect to~$\mu$ if 
	either of the two following equivalent conditions 
	is satisfied: 
	\begin{enumerate}[(i)]
		\item $\displaystyle
		\lim_{n\to \infty}
		\frac{1}{n}
		\sum_{k=0}^{n-1}
		| \mu ( A \cap T^{-k} (B)) - \mu (A) \mu (B) | = 0
		$,
		$A,B\in {\cal B}$,
		\item
		$\displaystyle
		\lim_{n\to \infty}
		\frac{1}{n}
		\sum_{k=0}^{n-1}
		\left|I_n(f,g)\right| =0$,
		$f,g \in L^2( E , \mu)$. 
	\end{enumerate} 
	When $\mu$ is a Gaussian measure on $E$, the mixing of linear operators $T$
	has been characterized in \cite[Theorem 5.24]{DLO} and references therein
	using the covariance operator $R:E^*\to E$ of $\mu$, defined as 
	\begin{equation}
		\nonumber 
		\langle Rx^*,y^*\rangle =
		\int_E \widebar{\langle z , x^* \rangle} \langle z, y^* \rangle \mu (dz),
		\qquad  x^*,y^* \in E^*, 
	\end{equation}
	where $E^*$ is the continuous dual of $E$ 
	and $\biganglebracket{\cdot, \cdot}: E \times E^*\to \C$
	denotes the dual product.  
		Such characterizations
	have been recently extended in \cite{MP23} from Gaussian measures to
	a wide class of infinitely divisible probability measures $\mu$ on
	{real and complex separable Banach spaces} $E$
	using strong and weak mixing properties
	of stationary infinitely divisible processes
	established in \cite{maruyama,SCMIDP,EEWMIDP, MCMIDP, MPMIDRF}. 
	
	\medskip
	
	Recall that a probability measure 
	$\mu$ on the Banach space~$E$
	is infinitely divisible if and only if 
	for every $n\geq 1$ there exists another probability measure $\mu_n$
	on $E$ such that 
	$$
	\mu =
	\underset{ n \mbox{\scriptsize{ times }}}{\underbrace{
			\mu_n \star \cdots \star \mu_n}}, 
	$$ 
	see e.g. \S5.1 of \cite{lindebook}, where $\star$ denotes
        measure convolution.
		It is known in addition that 
	every infinitely divisible probability measure on
	a complex Banach space~$E$
	has a characteristic functional of the form  
	\begin{align} 
	\label{jlkdf23} 	
 & 			
			\int_E e^{i \re \langle z , x^* \rangle} \mu (dz)
		\\
		\nonumber
		& \quad =  \exp \left(
		- \frac{1}{4} \biganglebracket{Rx^*,x^*}
		+ i \re \langle \hat{x} , x^* \rangle 
		+ \int_E \big(
		e^{i \re \biganglebracket{z,x^*}}  - 1 - i
		\kappa (z) \re \biganglebracket{z,x^*}
		\big)
		\lambda (dz) \right) \qquad  
	\end{align} 
	for all $x^*\in E^*$, 
	see e.g. \cite[\S II.1]{Rosinskialt}, where
	$\hat{x} \in E$, and 
	\begin{itemize}
		\item
		$R:E^*\to E$ is a conjugate symmetric and positive semidefinite
		covariance operator,
		\item $\lambda$ is a L\'evy measure, i.e.
		$\lambda$ is a measure on $E$ that satisfies
		$\lambda (\{0\}) = 0$ and 
                  \begin{equation}
                    \label{jkld145} 
		\int_E \min \big( 1 , (\re\biganglebracket{z,x^*})^2 \big) \lambda (dz) < \infty, \qquad x^*\in E^*,  
		\end{equation}
	      \item $\kappa (z)$ is a bounded measurable function on $E$ 
			such that $\lim_{z\to 0} \kappa (z) =1$ and
			$\kappa (z) = O(1/\bignorm{z})$ as $\bignorm{z}$ tends to infinity,
			called a truncation function.
	\end{itemize}
		In this infinitely divisible setting, the characterization result of \cite{MP23}
	uses codifference functionals 
	defined as
	\begin{equation}
		\nonumber 
		C^=_\mu (x^*, y^*) := \log
		\int_E
		{e^{i \re \biganglebracket{z,x^*-y^*}}}
		\mu (dz)
		- \log
		\int_E
		{e^{i \re \biganglebracket{z,x^*}}}\mu (dz)
		- \log
		\int_E
		{e^{-i \re \biganglebracket{z,y^*}}}
		\mu (dz)
		, 
	\end{equation} 
	and
	\begin{equation}
		\nonumber 
		C^{\neq}_\mu (x^*, y^*) :=
		\log \int_E
		{e^{
				i \re \biganglebracket{z,x^*}
				-
				i \imag  \biganglebracket{z,y^*}}
		}
		\mu (dz)
		- \log \int_E
		{e^{i \re \biganglebracket{z,x^*}}}
		\mu (dz)
		- \log \int_E
		{e^{-i \imag  \biganglebracket{z,y^*}}}
		\mu (dz) , 
	\end{equation} 
	$x^*,y^*\in E^*$. 
	
	\medskip
	
	In Theorem~\ref{fkds3} below we start by improving on 
	Proposition~2.2 of \cite{MP23}, 
	by removing the vanishing support Condition~\ref{fjklf34} 
	imposed therein on the L\'evy measure of the
	pushforwards of $\mu$ by linear functionals $x^*\in E^*$.
        This condition 
	originated in \cite{maruyama} and \cite{SCMIDP}, 
        and we rely on results of 
        \cite{MCMIDP}
        and
        \cite{MPMIDRF}
        who relaxed it in the framework of         
        stochastic processes. 
	As a result, we characterize the mixing
	of linear operators $T$ 
		via the asymptotic vanishing of codifferences 
	\begin{align}
		\label{a1}
		\lim\limits_{n\to\infty} C_\mu^=(ax^*, aT^{*n}x^*) = 0 \qquad\text{and} \qquad \lim\limits_{n\to\infty} C_\mu^{\neq}(ax^*, aT^{*n}x^*) = 0
	\end{align}
	for all $x^*$ and for some $a\not= 0$
	depending on $x^*\in E^*$,
	see Theorem~\ref{fkds3}. 
	In Examples~\ref{jkldf1}
	and \ref{poicor} 
	we consider measures $\mu$ that can be
	treated by
	Theorem~\ref{fkds3},
	and to which
	Proposition~2.2 of \cite{MP23} does not apply. 
	
	\medskip
	
	Next, we focus our attention on the speed of mixing via
	the derivation of decay rates for the codifferences appearing in 
	\eqref{a1}. 
	In the setting of Gaussian measures on Hilbert
	spaces $H$, covariance decay rates of the form
	\begin{align*}
		\bigabs{\biganglebracket{Rx^*, T^{*n}x^*}} \leq C 
		\frac{\|x^*\|^2}{n^\gamma}, 
	\end{align*}
	where $C$ is a constant independent of $x^*$,
	have been obtained in \cite{Devinck},
	provided that 
	$T$ is $\sigma$-spanning (i.e. for every $\sigma$-measurable subset $A\subset\T$ such that $\sigma(A)=1$, the eigenspaces $\ker(T-\lambda I), \lambda\in A$ span a dense subset of $H$),
	and $T$ admits a $\gamma$-H\"olderian
	eigenvector field for some $\gamma\in(0,1]$,
	where $\sigma$ denotes the normalized Lebesgue measure on 
	the complex unit circle. 
 In the more general setting of Banach spaces, similar covariance decay rates for
	classes of functions which satisfy a central limit theorem have been described in \cite{CLT}.
		
	\medskip
	
	In Section~\ref{s4} 
		we derive bounds on
	the codifferences $C_\mu^= (x^*, T^{*n} x^*)$
	and
	$C_\mu^{\neq} (x^*, T^{*n} x^*)$
	which provide quantitative estimates of 
	mixing speed 
	in the infinitely divisible setting. 
	For this, in addition to \eqref{jlkdf23}, we consider
	infinitely divisible measures with 
	characteristic functionals of the form 
	\begin{eqnarray}
		\label{fjkldf34} 
		\lefteqn{ 
			\int_E e^{i \re \langle z , x^* \rangle} \mu (dz)
		}
		\\
		\nonumber 
		& \quad = \exp \left( 
		\displaystyle - \frac{1}{4}
		\biganglebracket{Rx^*,x^*}
		+ 
		\int_E 
		\int_{-\infty}^\infty \big( e^{i u \re \langle z , x^* \rangle}-1-i u\kappa (u) \re\langle z , x^* \rangle
		\big) \rho(z , du) \xi (dz) \right) 
	\end{eqnarray} 
	for all $x^*\in E^*$, where
	\begin{itemize}
		\item $\{\rho( z , \cdot )\}_{z\in E}$ is
		a family of L\'evy measures on $\R$,
		\item $\xi $ is a $\sigma$-finite  measure on $E$ 
		called a control measure,
		and 
		\item $\kappa$ is the truncation function
		$$
		\kappa (u) := {\bf 1}_{\{ |u| < 1 \}}
		+ \frac{1}{|u|}
		{\bf 1}_{\{|u|\geq 1\}}, \quad
		u \in \real.
		$$
	\end{itemize}
        
	 In Section~\ref{s4.1},
        using the control measure bounds of Section~\ref{s4} 
	we derive mixing rates under compound Poisson measures 
	on $E=\ell^p (\N )$, $p\in[1,2)$, which have 
	characteristic functional 
	\eqref{jlkdf23} 
	and L\'evy measure of the form 
$$ 
		\displaystyle
		\lambda (dz) := \sum_{n=0}^\infty \delta_{\lambda_n e_n}(dz), 
                $$
 where $(e_n)_{n\geq 0}$ denotes the canonical basis of $\ell^p (\N )$
 and $\delta_x$ is the Dirac measure at any point $x \in E$.
 
 \medskip

 In Section~\ref{s4.2}
	we let 
	$\mu$ be an $\alpha$-stable measure 
	with $\alpha\in (0,2) \setminus \{1\}$,  
	in which case
	$\biganglebracket{Rx^*, y^*} = 0$ and 
	the characteristic functional 
	of $\mu$ 
		can be written
	by the Tortrat Theorem \cite{Tortrat} as
	\begin{align} 
		\label{st}
				\int_E e^{i\re  \langle z , x^* \rangle} \mu (dz)
		& = \exp\bigbracket{
			c_\alpha
			\int_E 
			\int_{-\infty}^\infty \big( e^{i u \re \langle z , x^* \rangle}-1-i u \kappa (u) \re\langle z , x^* \rangle
			\big) \frac{du}{|u|^{1+\alpha}} \xi (dz)}
		\\
		\nonumber
		& 
		= \exp\bigbracket{- \int_E \bigabs{\re  \biganglebracket{z,x^*}}^\alpha \xi (dz)}, \quad x^*\in E^*, 
	\end{align} 
	where $\xi$ is a finite
	control measure, 
	$\kappa (u)$ is a truncation function,
	and $c_\alpha  >0$, 
	see also \cite[Sec. 3]{SCMIDP}, \cite[p.~6]{GMBS}, \cite[Corollary~5.5]{LT}, 
	\cite[Theorem 6.4.4 
	and Corollary~7.5.2]{lindebook}, 
	\cite[Lemma~14.11]{sato}, or \cite[Corollary~4.1]{Pitmanarticle}. 

\medskip

 In particular, in Proposition~\ref{extexp} we derive explicit
codifference decay rates
of the form 
$$
\sup_{x^*,y^*\in E^*\setminus \{0\}}
\frac{
	\bigabs{C_\mu^{=,\neq} (x^*, T^{*n}y^*)}
}{\|x^*\|^{\alpha /2}\|y^*\|^{\alpha /2}
}
=
O\big(
\eta^{\alpha n/2}\big), 
$$ 
where
$C_\mu^{=,\neq}$ denotes 
$C_\mu^{=}$ or 
$C_\mu^{\neq}$, 
for weighted forward shifts $T$
on $E = \ell^p(\Z)$ 
under $\alpha$-stable measures 
with
$\alpha \in (1/2 , 2) \setminus \{1 \}$ and $p\in (\alpha , 2\alpha) \cap [1,2]$, where $\eta \in (0,1)$
  is a constant depending on $T$. 

  \medskip
   
 In Corollary~\ref{infseriesrate}, 
 we extend those results by deriving decay rates for the quantity
\begin{equation}
	\label{fjkl1} 
	I_n(f,g) :=
	\int_E f(z)g(T^n z) \mu(dz) - \int_E f(z)\mu(dz) \int_E g(z)\mu(dz), 
\end{equation}
where $f,g$ are finite or infinite linear
combinations of exponentials,
and $T$ is a weighted forward shift as in 
 Proposition~\ref{extexp}. 

\medskip

In Section~\ref{s6}
we consider the case where $\mu$ is a tempered stable measure
whose characteristic functional takes the form
\eqref{fjkldf34} 
and $\rho (z , du)$ is given by 
\begin{align*}
	\rho(z,du) =
	\bigbracket{
		\frac{a_-}{|u|^{1+\alpha }}e^{-\lambda_- |u|} \mathbf{1}_{\R^-}(u)
		+
		\frac{a_+}{u^{1+\alpha }}e^{-\lambda_+ u} \mathbf{1}_{\R^+}(u)}\du ,
\end{align*}
with
$a_-,a_+,\lambda_-,\lambda_+>0$ 
and $\alpha \in (0,1)$,
see \cite{TSDP} and also \cite{Devinck}.         
 In this setting we derive decay rates of the form 
$$ 
\sup_{x^*,y^*\in E^*\setminus \{0\}}
			\frac{
				\bigabs{C_\mu^{=,\neq} (x^*, T^{*n}y^*)}}
			     {\|x^*\|^{p/2}\|y^*\|^{p/2}}
                             = O ( n^{-\zeta }),
                             $$
                             for
                             certain backward weighted shift operators
          on 
                             $E =  \ell^p ( \Z )$, 
		             where
                             $\zeta \in (0,1)$ is a constant
                             depending on $T$
                             for codifferences and quantities of the form \eqref{fjkl1}, see Proposition~\ref{exampprop}
 and
 Example~\ref{exampcor}. 
\section{Mixing conditions} 
\noindent
The goal of this section is to prove Theorem~\ref{fkds3} below, which 
extends
the necessary and sufficient conditions for the
mixing of linear operators in terms of codifference functionals
in Proposition~2.2 of \cite{MP23} 
by removing the technical support Condition~\ref{fjklf34} below.

\medskip

Recall, see e.g. \cite[Theorem~1.2.14]{applebk2}, that similarly to \eqref{jlkdf23},
every infinitely divisible random variable on $\R^d$, $d\geq 1$, 
has a characteristic functional of the form  
\begin{equation} 
\nonumber 
\int_{\R^d} e^{i \langle z , y \rangle_d } \mu (dz)
=  \exp \left(
- \frac{1}{2} \biganglebracket{Ry,y}_d
+
i \langle y_0 , y \rangle_d 
+ \int_{\R^d} \big(
e^{i \biganglebracket{z,y}_d}  - 1 - i
\kappa (z) \biganglebracket{z,y}_d
\big)
\nu (dz) \right)
\quad 
\end{equation} 
 for all $y\in {\R^d}$, where 
\begin{itemize}
\item $y_0\in \R^d$ is a fixed vector,
\item $\langle \cdot , \cdot \rangle_d$ denotes the Euclidean
  inner product in $\real^d$, 
\item
$R:{\R^d}\to {\R^d}$ is a symmetric and positive semidefinite
covariance operator,
\item $\kappa (z)$ is a bounded measurable function on ${\R^d}$ 
such that $\lim_{z\to 0} \kappa (z) =1$ and
$\kappa (z) = O(1/\bignorm{z})$ as $\bignorm{z}$ tends to infinity,
called a truncation function, and
\item $\nu$ is a L\'evy measure on ${\R^d}$,
i.e. $\nu$ is a measure that satisfies
$\nu (\{0\}) = 0$ and 
$$
\int_{\R^d} \min \big( 1 , \biganglebracket{z,y}^2_d \big) \nu (dz) < \infty, \qquad y\in {\R^d}.
$$
\end{itemize}
 In what follows, we consider an $E$-valued random variable
 $X$ with infinitely divisible distribution $\mu$.
\begin{condition}
\label{fjklf34-0}
 For any $x^*\in E^*$, the L\'evy measure $\nu_{x^*}$ of
 the $\real^2$-valued random variable
 $(\re  \biganglebracket{X,x^*} , \imag  \biganglebracket{X,x^*})$ 
 satisfies 
\begin{equation}
\label{fjklf34} 
\nu_{x^*} ( \R \times 2\pi\Z ) = 0
\quad
and
\quad \nu_{x^*} ( 2\pi\Z \times \R  ) = 0. 
\end{equation}
\end{condition}
The above assumption originates from a condition appearing in \cite{maruyama}, 
 which was used in \cite[Theorem~1]{SCMIDP}
to characterize mixing by codifferences.
The value of $2\pi$ in \eqref{fjklf34}
is chosen for consistency with the literature,
however, it can be replaced
with an arbitrary non-zero constant without affecting
 Condition~\ref{fjklf34}.
\begin{definition}
Given $\nu$ a measure on $\real^d$, 
we define the set $Z_d ( \nu )$ as follows.
If
\begin{equation}
\nonumber 
\nu ( \real^{j-1}\times \{2 k \pi \} \times \real^{d-j})=0 
\ \text{for all} \ k\in \Z, \ j = 1, \ldots , d, 
\end{equation}
we define $Z_d (\nu)  := \R \setminus \{1\}$,
else we let 
\begin{align*}
Z_d ( \nu ) := \bigcup_{k\in\Z}
\left\{k \frac{2 \pi}{s} \ : \ 
s\in \real\setminus \{0\} \ 
\mbox{and}
\ \sum_{j=1}^d
\nu ( \real^{j-1}\times \{s\} \times \real^{d-j}) > 0
\right\}. 
\end{align*}
\end{definition} 
\noindent
 Lemma~\ref{genscaling-0} shows that
for any L\'evy measure $\nu$, the set 
$\R\backslash Z_d ( \nu ) $ always contains a non-zero element.
\begin{lemma}
\label{genscaling-0}
For any L\'evy measure $\nu$  on $\real^d$, 
the set $Z_d ( \nu ) $ is either $\R\backslash\{1\}$ or is at most countable. 
\end{lemma}
\begin{Proof}
Clearly, we may assume that
\begin{align}
\nonumber 
\sum_{j=1}^d \nu ( \real^{j-1}\times \{2 k \pi \} \times \real^{d-j}) > 0  
\end{align} 
for some $k\in \Z$,
otherwise $Z_d ( \nu ) = \R\backslash\{1\}$ and
the proof is concluded. 
Denoting by $B_d(0,1)$ the unit ball and by
$\|\cdot\|_d$ the Euclidean norm of $\R^d$,
by definition of L\'evy measures we have 
\begin{align*}
\int_{\R^d} \min(1, \|x\|_d^2) \nu (dx) = \int_{B_d(0,1)} \|x\|_d^2 \nu (dx) + \int_{\R^d\backslash B_d(0,1)} \nu (dx) < \infty. 
\end{align*}
Letting
$P_1:\R^d\to \R$
denote the projection onto the first coordinate in
$\R^d$, 
and denoting by $\rho$ the pushforward of $\nu $ by $P_1$,
we have
\begin{align*}
&  \int_{-\infty}^\infty \min(1, x^2) \rho(dx) = \int_{(-1,1)} x^2 \rho(dx) + \int_{\R\backslash(-1,1)} \rho(dx)
\\
& \quad  \quad \quad
\leq 
\int_{(-1,1)\times\R^{d-1}} (P_1x)^2 \nu (dx)
+ \int_{\R^d\backslash B_d(0,1)} \nu (dx)
\\
& \quad \quad \quad= \int_{B_d(0,1)} (P_1x)^2 \nu (dx) + \int_{((-1,1)\times\R^{d-1}) \backslash B_d(0,1)} (P_1x)^2 \nu (dx)
+ \int_{\R^d\backslash B_d(0,1)} \nu (dx)
\\
&   \quad \quad\quad\leq \int_{B_d(0,1)} \|x\|_d^2 \nu (dx) + \int_{((-1,1)\times\R^{d-1}) \backslash B_d(0,1)} \nu (dx)
+ \int_{\R^d\backslash B_d(0,1)} \nu (dx)
\\
& \quad \quad \quad< \infty, 
\end{align*}
hence $\rho$ is $\sigma$-finite
and therefore it has countably many atoms,
i.e. there are at most countably many
those values of $s\in \real$ such that
$		  \nu (\{s\} \times \real^{d-1} ) > 0$.
Repeating this argument for each coordinate, we find that the cardinality of $Z_d ( \nu ) $ is at most $\N^d$, i.e. countable.
\end{Proof}
\noindent 
 The following results are stated for complex Banach spaces,
but they also apply to real Banach spaces  
by ignoring vanishing imaginary components.
Theorem~\ref{fkds3} 
allows for the mixing property of $T$ to be
checked without imposing Condition~\ref{fjklf34-0}.
Recall that a set $D\subset \N$ has density one if
$\lim_{n\to \infty} |D\cap \{0,1,\ldots , n\}|/(n+1) = 1$. 
\begin{theorem}
\label{fkds3}
Let $\mu$ be an infinitely divisible distribution
on a complex separable Banach space $E$.
For any $x^*\in E^*$, let $\nu_{x^*}$ denote  the L\'evy measure of
 $( \re \biganglebracket{X,x^*}, \imag \biganglebracket{X,x^*})$
on $\real^2$. 
Then, 
\begin{enumerate}[i)]
\item
$T$ is mixing if and only if for each $x^*\in E^*$ we have  
	\begin{align*}
		\lim\limits_{n\to\infty} C_\mu^=(ax^*, aT^{*n}x^*) = 0 \qquad\text{and} \qquad \lim\limits_{n\to\infty} C_\mu^{\neq}(ax^*, aT^{*n}x^*) = 0
	\end{align*}
	for some non-zero $a\in\R\backslash Z_2 ( \nu_{x^*})$. 
	\item $T$ is weakly mixing if and only if
for each $x^*\in E^*$
there exists a density one set $D_{x^*}\subset\N$ such that 
	\begin{align*}
		\lim\limits_{n\to\infty \atop n\in D_{x^*}} C_\mu^=(ax^*, aT^{*n}x^*) = 0 \qquad\text{and} \qquad \lim\limits_{n\to\infty \atop n\in D_{x^*}} C_\mu^{\neq}(ax^*, aT^{*n}x^*) = 0
	\end{align*}
	for some non-zero $a\in\R\backslash Z_2(\nu_{x^*})$.
	\end{enumerate}
\end{theorem}
In Example~\ref{jkldf1}
we present a non-mixing operator 
inspired by the example on 
page~282 of
\cite{SCMIDP},
which can be treated by Theorem~\ref{fkds3}
while Condition~\ref{fjklf34-0} is not satisfied. 
See also Example~\ref{poicor} for a mixing operator. 
\begin{example}
\label{jkldf1}
 Let $E$ be the (real) sequence space $\ell^p (\N )$, 
 $p\geq 1$, take 
$$
R:=0, \ 
\hat{x} := 2\pi e_0, \ 
\kappa \mbox{ such that } \kappa ( \hat{x} ) = 1,
 \mbox{ and let } \lambda (dz) := \delta_{2\pi e_0}(dz)
\mbox{ in } \eqref{jlkdf23}.
$$ 
Then, for $x^*\in E^*$ such that $\langle e_0,x^*\rangle = 1$, 
the L\'evy measure $\nu_{x^*} = \delta_{2\pi \langle e_0,x^*\rangle}$ on $\real$ 
does not satisfy \eqref{fjklf34}. 
However, the application of Theorem~\ref{fkds3}
 confirms that the identity operator $T={\rm Id}$ is 
not mixing. 
\end{example}
\begin{Proof} 
Let $\mu$ be the infinitely
divisible distribution on $E = \ell^p(\N)$ 
with L\'evy measure 
$\lambda (dz)$. 
In this case, $X$ can be defined as
 $X := 2\pi N e_0$ where $N$ is a standard Poisson random variable, and
 we have
 \begin{align*} 
 \E \big[ 
   {e^{i a \biganglebracket{X,x^*}}}
   \big] 
 & =  
\int_E 
{e^{i a \biganglebracket{z,x^*}}}\mu (dz)
\\
 & =  \exp \left(  \int_E \big(
e^{i a \biganglebracket{z,x^*}}  - 1 
\big)
\lambda (dz) 
\right)
\\
 & =  
\exp \left(
e^{2ia\pi \biganglebracket{e_0,x^*}}  - 1 \right), 
\end{align*} 
hence
for any $x^*\in E^*$
the random variable
 $\langle X,x^*\rangle = 2\pi N \langle e_0 , x^*\rangle$
 has L\'evy measure $\nu_{x^*} = \delta_{2\pi \langle e_0,x^*\rangle}$. 
In this case, we have
$$
Z_1(\nu_{x^*}) = 
\left\{ \frac{k}{\langle e_0,x^*\rangle} \ : \ k\in \Z\right\}
$$ 
if $\langle e_0,x^*\rangle \notin \Z$, and
$Z_1(\nu_{x^*}) = \real \setminus \{1\}$ otherwise. 
Hence, 
$$ 
C_\mu^=(ax^*, aT^{*n}x^*) =  
C_\mu^=(ax^*, ax^*)
=
-2 \log \int_E 
{e^{i a \biganglebracket{z,x^*}}}\mu (dz)
=
-2 \left(
e^{2ia\pi \biganglebracket{e_0,x^*}}  - 1 \right)
$$ 
is constant in $n\geq 1$
and does not vanish for any $a\in \real \setminus Z_1(\nu_{x^*})$, 
therefore Theorem~\ref{fkds3} confirms that $T$ is not mixing. 
\end{Proof} 
The proof of Theorem~\ref{fkds3} is 
stated at the end of this section
 by carrying over Theorem~2 of \cite{SCMIDP}
from the stochastic process setting 
to the framework of linear dynamics 
thereby completing the characterization of mixing of
infinitely divisible measures on Banach spaces. 
 For this, we need to prove the following multidimensional extension of
 \cite[Theorem 2]{SCMIDP}
 on mixing and weak mixing for discrete-time stochastic processes,
 which removes the support
condition assumed in \cite[Theorem 2.1]{MCMIDP}
and \cite[Theorem 4.3]{MPMIDRF}. 
\begin{prop}
\label{genscaling}
Let $d\geq 1$,
and let
$(X_n)_{n\geq 0}
= \big(X_n^{(1)},\ldots , X_n^{(d)}\big)_{n\geq 0}
$ 
be a stationary infinitely
divisible $\R^d$-valued process.
Denote by $\nu_0$ the L\'evy measure of $X_0$.
Then, 
\begin{enumerate}[a)]
\item
$(X_n)_{n\geq 0}$ is 
mixing if and only if for some non-zero $a\in\R\backslash Z_d ( \nu_0 )$ we have
\begin{align*}
	\lim\limits_{n\to\infty} \E \big[
	e^{ia(X_n^{(j)} - X_0^{(k)})}\big] = \E \big[
	e^{iaX_0^{(j)}}\big] \E \big[ e^{-iaX_0^{(k)}}
	\big],
\end{align*}
for any $j,k\in\{1,\dots, d\}$; and  
\item
$(X_n)_{n\geq 0}$ is weakly mixing if and only if for some non-zero $a\in\R\backslash Z_d (\nu_0 )$ and a density one set $D\subset\N$ we have
\begin{align*}
	\lim\limits_{n\to\infty \atop n\in D} \E
	\big[
	e^{ia(X_n^{(j)} - X_0^{(k)})
	}\big]
	= \E \big[
	e^{iaX_0^{(j)}}\big]
	\E \big[ e^{-iaX_0^{(k)}} \big],
\end{align*}
for any $j,k\in\{1,\dots, d\}$.
\end{enumerate}
\end{prop}
\begin{Proof}
If the condition 
\begin{align}\label{zerosupp}
\nu_0\bigbracket{\{x = (x_1,\dots, x_d)\in\R^d \ : \ \exists j\in \{1,\dots, d\}, \ x_j\in 2\pi\Z\}} = 0
\end{align}
holds, we can conclude from
\cite[Theorem 2.1]{MCMIDP} in the mixing case and
from \cite[Theorem 4.3]{MPMIDRF} in the weak mixing case,
hence we may assume that \eqref{zerosupp} does not hold.
Then, as in \cite[Theorem 2]{SCMIDP}
we observe that $(X_n)_{n\geq 0}$ is mixing, respectively weak mixing, if and only if $(aX_n)_{n\geq 0}$ is, and furthermore the L\'evy measure $\nu_0^a(\cdot)$ of $aX_0$ is $\nu_0(a^{-1}(\cdot))$. Now, since $a\notin Z_d ( \nu_0 )$, it follows that
\begin{align*}
\nu_0^a\bigbracket{\big\{x = (x_1,\dots, x_d)\in\R^d \ : \ \exists j\in \{1,\dots, d\} \mbox{ s.t. } x_j\in 2\pi\Z \big\}} = 0.
\end{align*}
The conclusion follows for mixing by \cite[Theorem 2.1]{MCMIDP}, see also \cite[Theorem 3.2]{MPMIDRF}, and for weak mixing by \cite[Theorem 4.3]{MPMIDRF}.
\end{Proof}
\noindent 
\begin{Proofy}{\em of Theorem~\ref{fkds3}.} 
Let $X$ denote a random variable with distribution $\mu$ on $E$.
For any $x^*\in E^*$, let the process
$(X^{x^*}_n)_{n\geq 0}$ be defined by 
\begin{equation} 
\nonumber 
X^{x^*}_n 
:= \big( \re \langle X, T^{*n}x^* \rangle ,
\imag  \langle X, T^{*n}x^* \rangle \big), 
\qquad n\geq 0. 
\end{equation}
By \cite[Lemma 2.1]{MP23}, $T$ is mixing if and only if $(X_n^{x^*})_{n\geq 0}$ is mixing for each $x^*\in E^*$, 
hence we conclude from Proposition~\ref{genscaling}
and the relations
\begin{align*} 
&
\qquad \qquad
C_{\mu}^=(ax^*, aT^{*n}x^*) = \log
\frac{
  \E\bigsquarebracket{
    e^{
      ia\re\biganglebracket{X, x^* } 
- 
      ia\re\biganglebracket{X, T^{*n}x^*
    }}
  }
}{ 
	\E\bigsquarebracket{e^{ia\re\biganglebracket{X, x^*}}}
	\E\bigsquarebracket{e^{-ia\re\biganglebracket{X, x^*}}}
} \qquad \qquad \qquad \qquad \qquad \qquad \qquad \qquad 
\\
& \mbox{and}  
\\
&  
\qquad \qquad  C_{\mu}^{\neq}(ax^*, aT^{*n}x^*)
= \log
\frac{
  \E\bigsquarebracket{e^{
      ia\re\biganglebracket{X, x^*}
      - ia\imag\biganglebracket{X, T^{*n}x^*}}}}{
	\E\bigsquarebracket{e^{ia\re\biganglebracket{X, x^*}}}
	\E\bigsquarebracket{e^{-ia\imag\biganglebracket{X, x^*}}}
}, \qquad
a\in\R. 
\end{align*} 
\end{Proofy}
\section{Codifference bounds} 
\label{s4}
\noindent 
Our first step towards the derivation of
codifference decay rates in \eqref{a1} is to
derive bounds on codifferences
using L\'evy and control measures. 
\noindent
\begin{lemma}[L\'evy measure bounds.]
\label{lemma3.1}
Let $\mu$ be an infinitely divisible distribution with
characteristic functional of the form \eqref{jlkdf23}. 
For every $p\in [0,2]$, we have 
the codifference bounds
\begin{align}
	\label{b2-00}
	\bigabs{C_\mu^= (x^*, y^*)} \leq \frac{1}{2}\bigabs{\re \biganglebracket{Rx^*, y^*}} + 2^{4-p} \int_E |\re\biganglebracket{z,x^*}\re\biganglebracket{z,y^*}|^{p/2}\lambda(dz),
\end{align}
and
\begin{align}
	\label{b2}
	\bigabs{C_\mu^{\neq}(x^*, y^*)} \leq \frac{1}{2}\bigabs{\imag \biganglebracket{Rx^*, y^*}} + 2^{4-p} \int_E |\re\biganglebracket{z,x^*}\imag\biganglebracket{z,y^*}|^{p/2}\lambda(dz).
\end{align}
\end{lemma}
\begin{Proof}
\noindent
 The codifference of $\mu$ can be rewritten from \eqref{jlkdf23} as 
\begin{equation}
	\label{cdf1} 
	C^=_\mu (x^*, y^*)
	=  
	\frac{1}{2}\re \biganglebracket{Rx^*,y^*} + \int_E
	\big(
	e^{i\re \biganglebracket{z,x^*}}  - 1 
	\big)
	\big(
	e^{ - i \re \biganglebracket{z,y^*}}  - 1 
	\big)
	\lambda (dz), 
\end{equation} 
and
\begin{equation}
	\label{cdf2} 
	C^{\neq}_\mu (x^*, y^*)
	=  
	\frac{1}{2}\imag \biganglebracket{Rx^*,y^*} + \int_E
	\big(
	e^{i\re \biganglebracket{z,x^*}}  - 1 
	\big)
	\big(
	e^{ - i \imag  \biganglebracket{z,y^*}}  - 1 
	\big)
	\lambda (dz), 
\end{equation} 
$x^*,y^*\in E^*$. 
Taking the real part in \eqref{cdf1}, we have
\begin{align*}
	&	\re C_\mu^= (x^*, y^*) 
	= \frac{1}{2}\re\biganglebracket{Rx^*, y^*}\\
	&\quad + \int_E
	\big(
	(\cos(\re\biganglebracket{z,x^*}) - 1)(\cos(\re\biganglebracket{z,y^*}) - 1) + \sin(\re\biganglebracket{z,x^*})\sin(\re\biganglebracket{z,y^*})
	\big)
	\lambda(dz).
\end{align*}
Likewise, taking the imaginary part in \eqref{cdf1}, we obtain
\begin{align*}
	& \imag C_\mu^= (x^*, y^*)
	\\
	& \qquad
	= \int_E
	\big(
	\sin(\re\biganglebracket{z,x^*})(\cos(\re\biganglebracket{z,y^*}) - 1) - \sin(\re\biganglebracket{z,y^*})(\cos(\re\biganglebracket{z,x^*}) - 1)
	\big)\lambda(dz).
\end{align*}
Using the inequalities
\begin{align}
	\label{fjklsd241} 
	\max ( |\cos x - 1|, |\sin x| ) \leq
	2^{\frac{2-p}{2}} |x|^{p/2}, \quad x\in \real, 
\end{align}
which are valid for $p\in [0,2]$,
it follows from the triangle inequality that
\begin{align*}
	\bigabs{\re C_\mu^= (x^*, y^*)} \leq \frac{1}{2}\bigabs{\re \biganglebracket{Rx^*, y^*}} + 2^{3-p}
	\int_E |\re\biganglebracket{z,x^*}\re\biganglebracket{z,y^*}|^{p/2} \lambda(dz).
\end{align*}
Similarly, we obtain
\begin{align*}
	\bigabs{\imag C_\mu^= (x^*, y^*)} \leq
	2^{3-p} \int_E |\re\biganglebracket{z,x^*}\re\biganglebracket{z,y^*}|\lambda(dz), 
\end{align*}
which yields \eqref{b2-00}. 
Likewise, for \eqref{cdf2} we have 
\begin{align*}
	& \re C_\mu^{\neq}(x^*, y^*)
	= \frac{1}{2}\re\biganglebracket{Rx^*, y^*}
	\\
	& \quad + \int_E
	\big(
	\cos(\re\biganglebracket{z,x^*}) - 1)(\cos(\imag\biganglebracket{z,y^*}) - 1) + \sin(\re\biganglebracket{z,x^*})\sin(\imag\biganglebracket{z,y^*}
	\big)\lambda(dz), 
\end{align*}
implying 
\begin{align*}
	\bigabs{\re C_\mu^{\neq}(x^*, y^*)} \leq \frac{1}{2}\bigabs{\imag \biganglebracket{Rx^*, y^*}} + 2^{3-p} \int_E |\re\biganglebracket{z,x^*}\imag\biganglebracket{z,y^*}|^{p/2} \lambda(dz), 
\end{align*}
and
\begin{align*}
	& \imag C_\mu^{\neq}(x^*, y^*)
	\\
	&\quad = \int_E
	\big(
	\sin(\re\biganglebracket{z,x^*})(\cos(\imag\biganglebracket{z,y^*}) - 1) - \sin(\imag\biganglebracket{z,y^*})(\cos(\re\biganglebracket{z,x^*}) - 1)
	\big)
	\lambda(dz), 
\end{align*}
implying 
\begin{align*}
	\bigabs{\imag C_\mu^{\neq}(x^*, y^*)} \leq
	2^{3-p}
	\int_E |\re\biganglebracket{z,x^*}\imag\biganglebracket{z,y^*}|^{p/2}
	\lambda(dz), 
\end{align*}
which yields \eqref{b2}.
\end{Proof}
\noindent
Next, we consider the case where 
the characteristic functional  of
the infinitely divisible measure
$\mu$ takes the form \eqref{fjkldf34}.
Recall, see \S3 of \cite{SCMIDP}, that any random variable $X$ 
with distribution $\mu$ can be represented as
$$
X= \int_E z\Lambda(dz)
$$
where 
$\Lambda$ is the infinitely divisible random measure
on $E$ defined by its characteristic functional  
$$ 
\E\bigsquarebracket{e^{i t \Lambda(A)}} \\
= \exp\bigbracket{
- \frac{
	t^2}{4} \int_A \sigma^2(z)\xi(dz)
+ \int_A \int_{-\infty}^\infty \big( e^{i u t } - 1 - i t u \kappa (u)
\big) \rho(z , du) \xi (dz)}
$$ 
for measurable $A\subset E$ and $t \in \R$, 
where $\sigma^2:E\to [0,\infty )$ is a measurable function. 
\begin{lemma}[Control measure bounds.] 
\label{fkjdf43}
Let $\mu$ be an infinitely divisible distribution with
characteristic functional of the form \eqref{fjkldf34}.
For any $p\in [0,2]$ and $c>0$, we have the codifference bounds 
\begin{align}
	\label{e1}
	\bigabs{C_\mu^= (x^*, y^*)}
	& \leq \frac{1}{2}\bigabs{\re \biganglebracket{Rx^*, y^*}}
	\\
	\nonumber
	& \quad + 16\int_E \bigbracket{2^{-p} \int_{-c}^c |u|^p |\re\biganglebracket{z,x^*}\re\biganglebracket{z,y^*}|^{p/2} \rho(z,du) +
		\rho(z,{\R\backslash[-c, c]})}\xi(dz).
\end{align}
and 
\begin{align}
	\label{e2}
	\bigabs{C_\mu^{\neq}(x^*, y^*)} & \leq \frac{1}{2}\bigabs{\imag \biganglebracket{Rx^*, y^*}}
	\\
	\nonumber
	& \quad + 16\int_E \bigbracket{2^{-p} \int_{-c}^c |u|^p |\re\biganglebracket{z,x^*}\imag\biganglebracket{z,y^*}|^{p/2} \rho(z,du) + \rho(z,{\R\backslash[-c, c]})}\xi(dz), 
\end{align}
for    $x^*,y^*\in E^*$.
\end{lemma}
\begin{Proof} 
By taking the real part
in the relation 
\begin{equation}
	\nonumber 
	C^=_\mu (x^*, y^*)
	=  
	\frac{1}{2}\re \biganglebracket{Rx^*,y^*} + \int_E \int_{-\infty}^\infty 
	\big(
	e^{i u \re \biganglebracket{z,x^*}}  - 1 
	\big)
	\big(
	e^{ - i u \re \biganglebracket{z,y^*}}  - 1 
	\big)\rho(z , du) \xi (dz), 
\end{equation} 
$x^*,y^*\in E^*$, we have
\begin{align*}
	& \re C_\mu^= (x^*, y^*) = \frac{1}{2}\re\biganglebracket{Rx^*, y^*}
	\\
	& 
	+
	\hskip-0.1cm
	\int_E \int_{-\infty}^\infty
	\hskip-0.2cm
	\big(
	(\cos(u\re\biganglebracket{z,x^*}) - 1 )(\cos(u\re\biganglebracket{z,y^*}) - 1) + \sin(u\re\biganglebracket{z,x^*})\sin(u\re\biganglebracket{z,y^*})
	\big)\rho(z,du)\xi(dz).
\end{align*}
 For any $c > 0$, using \eqref{fjklsd241}, we have
\begin{align*}
	& \int_E \int_{-\infty}^\infty \bigabs{(\cos(u\re\biganglebracket{z,x^*}) - 1)(\cos(u\re\biganglebracket{z,y^*}) - 1)} \rho(z,du)\xi(dz)\\
	& \quad \leq \int_E \bigbracket{2^{2-p}\int_{-c}^c |u|^p |\re\biganglebracket{z,x^*}\re\biganglebracket{z,y^*}|^{p/2} \rho(z,du)
		+ 4 \int_{\R\backslash[-c, c]} \rho(z,du)}\xi(dz).
\end{align*}
Likewise, we have
\begin{align*}
	& \int_E \int_{-\infty}^\infty \bigabs{\sin(u\re\biganglebracket{z,x^*})\sin(u\re\biganglebracket{z,y^*})} \rho(z,du)\xi(dz)\\
	& \quad \leq \int_E \bigbracket{2^{2-p}\int_{-c}^c  |u|^p |\re\biganglebracket{z,x^*}\re\biganglebracket{z,y^*}|^{p/2} \rho(z,du)
		+ 4 \int_{\R\backslash[-c, c]} \rho(z,du)}\xi(dz),
\end{align*}
hence by the triangle inequality we have
\begin{align*}
	&  \bigabs{\re C_\mu^= (x^*, y^*)}  \leq \frac{1}{2}\bigabs{\re \biganglebracket{Rx^*, y^*}} \\
	& \qquad \qquad
	\qquad
	+ 8\int_E \bigbracket{2^{-p} \int_{-c}^c |u|^p |\re\biganglebracket{z,x^*}\re\biganglebracket{z,y^*}|^{p/2} \rho(z,du) + \int_{\R\backslash[-c, c]} \rho(z,du)}\xi(dz).
\end{align*}
In a similar fashion, we have
\begin{align*}
	\bigabs{\imag C_\mu^= (x^*, y^*)} \leq 8\int_E \bigbracket{2^{-p} \int_{-c}^c  |u|^p |\re\biganglebracket{z,x^*}\re\biganglebracket{z,y^*}|^{p/2} \rho(z,du) + \int_{\R\backslash[-c, c]} \rho(z,du)}\xi(dz),
\end{align*}
which yields \eqref{e1}.
The bound \eqref{e2}
is obtained by application of similar arguments to 
\begin{equation}
	\nonumber 
	C^{\neq}_\mu (x^*, y^*)
	=  
	\frac{1}{2}\imag \biganglebracket{Rx^*,y^*} + \int_E \int_{-\infty}^\infty 
	\big(
	e^{i u \re \biganglebracket{z,x^*}}  - 1 
	\big)
	\big(
	e^{ - i u \imag  \biganglebracket{z,y^*}}  - 1 
	\big)\rho(z , du) \xi (dz). 
\end{equation} 
\end{Proof}

\section{Compound Poisson measures}
\label{s4.1}
\noindent
 In this section, we provide an example of 
mixing operator $T$ in Proposition~\ref{poiprop}
that can be treated
by Theorem~\ref{fkds3}, and to which
Proposition~2.2 of \cite{MP23} does not apply,
see Example~\ref{poicor} below. 
For this, we consider an infinitely divisible measure 
$\mu$ on $E=\ell^p (\N )$
 with characteristic functional \eqref{jlkdf23}
 and L\'evy measure of the form
	\begin{align}\label{leeeeeevy}
		\sum_{n=0}^\infty \delta_{\lambda_n e_n}(dz),
	\end{align}
	for appropriate sequences $(\lambda_n)_{n\geq 0}$,
        i.e. $\mu$ is the distribution of
        an $E$-valued random variable
        $X$ whose components in the canonical basis
        $(e_n)_{n\geq 0}$ of $\ell^p (\N )$
        are independent Poisson random variables
        with means $(\lambda_n)_{n\geq 0}$. 
	We first determine precisely the sequences $(\lambda_n)_{n\geq 0}$ for which \eqref{leeeeeevy} defines a L\'evy measure on $\ell^p(\N)$, as a consequence of the following auxiliary lemma.
	\begin{lemma}\label{holdup}
	  Let $p\in[1,2)$ and $q$ be the H\"older conjugate of $p$.
          For any complex sequence $(a_n)_{n\geq 0}$,
           the following statements are equivalent.
		\begin{enumerate}
			\item $(a_n)_{n\geq 0} \in \ell^p(\N)$.
			\item For every $(b_n)_{n\geq 0}$ in the real sequence
                           space $\ell^q(\N)$, we have
			  \begin{align*}
				\sum_{n=0}^\infty |a_n|^2 | b_n|^2 < \infty.
			\end{align*}
		\end{enumerate}
	\end{lemma}
	\begin{Proof}
	 The implication (1~$\Rightarrow$~2) follows from
                    H\"older's inequality 
		\begin{align*}
		  \sum_{n=0}^\infty |a_n| |b_n|
                                    \leq \|(a_n)_{n\geq 0}\|_{\ell^p(\N)} \|(b_n)_{n\geq 0}\|_{\ell^q(\N)} < \infty
		\end{align*}
	        and the fact that 
                $\ell^1(\N) \subset \ell^2(\N)$.
		The implication (2~$\Rightarrow$~1) follows from
                                the reverse H\"older inequality 
		\begin{align*}
			\infty > \sum_{n=0}^\infty |a_n|^2 |b_n|^2 \geq \bigbracket{\sum_{n=0}^\infty |a_n|^{2/s}}^s \bigbracket{\sum_{n=0}^\infty |b_n|^{-2/(s-1)}}^{-(s-1)}
		\end{align*}
                applied with $s:=2/p\in [1,2)$ 
                to any
                sequence
                $(b_n)_{n\geq 0} \in \ell^q(\N )$
                such that $\re(b_n) \neq 0$ for all $n\geq 0$. 
	\end{Proof}
	From Lemma~\ref{holdup}, we obtain the following criterion for L\'evy measures in real $\ell^p(\N)$.
	\begin{lemma}\label{propsy}
	  Let $p\in [1,2)$, and let $(\lambda_n)_{n\geq 0}$
          be a real sequence. The measure $\lambda$ defined by
		\begin{align*}
			\lambda (dz) := \sum_{n=0}^\infty \delta_{\lambda_n e_n}(dz)
		\end{align*}
		is a L\'evy measure on the real sequence space
                 $\ell^p(\N)$ if and only if $(\lambda_n)_{n\geq 0} \in \ell^p(\N)$ and $\lambda_n \neq 0$ for all $n\geq 0$.
	\end{lemma}
	\begin{Proof}
	   Let $q>2$ denote the H\"older conjugate of $p$. We note that 
		\begin{align*}
			\int_{\ell^p(\N)} \min(1,\biganglebracket{z,x^*}^2) \lambda(dz) = \sum_{n=0}^\infty \min(1,\lambda_n^2\biganglebracket{e_n,x^*}^2)
		\end{align*}
 is finite for all $x^* \in \ell^q (\N)$ if and only if 
		\begin{align*}
                                    \sum_{n=0}^\infty \lambda_n^2\biganglebracket{e_n,x^*}^2
 		\end{align*}
 is finite for any $x^* \in \ell^q (\N)$, i.e. if and only if 
                		\begin{align*}
                                    \sum_{n=0}^\infty \lambda_n^2 b_n^2 < \infty
		\end{align*}
                                for any $(b_n)_{n\geq 0}$
                                in the real sequence
                           space $\ell^q(\N )$.
                                We conclude
                                according to \eqref{jkld145}, using Lemma~\ref{holdup} and
		                the fact that $\lambda(\{0\}) = 0$
                                if and only if $\lambda_n \neq 0$ for all $n\geq 0$. 
                                \end{Proof} 
          \noindent
          We now turn to the main result of this section. 
\begin{prop}
	\label{poiprop}
	Let $p\in [1,2)$.
        Let $(\omega_n)_{n\geq 0}$ be a positive weight sequence
	such that the sequence $(\lambda_n)_{n\geq 0}$ defined by $\lambda_0 > 0$ and
	\begin{align*}
		\lambda_n: = \lambda_0 \prod_{l=0}^{n-1}
		\frac{1}{\omega_l},
                \quad 
		n\geq 0, 
	\end{align*}
	satisfies $(\lambda_n)_{n\geq 0}\in\ell^p(\N)$.
	Let $\mu$ be the compound Poisson measure on the real space
        $E: = \ell^p(\N)$ with characteristic functional 
	\eqref{jlkdf23} given by 
	$$
        R:=0,
        \quad 
	\hat{x} :=
	\sum_{n=0}^\infty
	\lambda_n e_n 
        \ \ 
        \mbox{and}
        \ \ \kappa (z) : = {\bf 1}_{\big\{
          \| z\|_{\ell^p (\N)} \leq \max\limits_{n\geq 0} \lambda_n\big\}}, \ z\in \ell^p (\N), 
                $$
        and the L\'evy measure
	\begin{equation}
		\label{fjk2}
		\displaystyle
		\lambda (dz) := \sum_{n=0}^\infty \delta_{\lambda_n e_n}(dz). 
	\end{equation}	
	Consider the weighted backward shift operator 
	$T$              defined by 
	$Te_0 := 0$ and 
	$$Te_{n+1} := \omega_n e_n, \quad n\geq 0.
	$$
	The following are true.
        \begin{enumerate}[1.]
        \item $T$ admits $\mu$ as invariant measure.
        \item
           We have 
	$$
	\sup_{x^*, y^*\in E^*\setminus \{0\}}
	\frac{\bigabs{C_\mu^= (x^*, T^{*n}y^*)}
	}{  \|x^*\|^{p/2}\|y^*\|^{p/2}}
	\leq
	2^{4-p}
	\sum_{l=0}^\infty
	( \lambda_l \lambda_{l+n}  )^{p/2},
	\quad x^*,y^*\in E^*,
	\ n \geq 0.
	$$ 
	 In particular, $T$ is mixing
	by Theorem~\ref{fkds3},             provided that
	$\displaystyle \lim\limits_{n\to\infty}
	\sum_{l=0}^\infty 
	( \lambda_l \lambda_{l+n} )^{p/2} = 0$.
        \end{enumerate}
\end{prop}
\begin{Proof}
  \noindent
   1. Note that \eqref{fjk2} defines a L\'evy measure by Lemma~\ref{propsy}, thus $\mu$ is well-defined.
	To show that $T$ admits $\mu$ as invariant measure, we
	use \eqref{jlkdf23} and the equalities 
		\begin{align*} 
		&\int_E \big(
		e^{i \re \biganglebracket{Tz,x^*}}  - 1  - i
		\kappa (z) \re \biganglebracket{Tz,x^*} \big)
		\lambda (dz)\\
		& \quad =
		\sum_{n=0}^\infty
		\int_E \big(
		e^{i \re \biganglebracket{Tz,x^*}}  - 1 - i
		\kappa (z) \re \biganglebracket{Tz,x^*}   \big)
		\delta_{\lambda_n e_n}(dz)
		\\
		& \quad =
		\sum_{n=1}^\infty
		\big(
		e^{i \re \biganglebracket{\lambda_n Te_n,x^*}}  - 1 - i
		\re \biganglebracket{\lambda_n Te_n,x^*}   \big)
		\\
		& \quad =
		\sum_{n=0}^\infty
		\big(
		e^{i \re \biganglebracket{\lambda_{n+1} Te_{n+1},x^*}}  - 1  - i
		\re \biganglebracket{\lambda_{n+1} Te_{n+1},x^*}  \big)
		\\
		& \quad =
		\sum_{n=0}^\infty
		\big(
		e^{i \re \biganglebracket{\lambda_n e_n,x^*}}  - 1  - i
		\re \biganglebracket{\lambda_n e_n,x^*} \big)
		\\
		& \quad =   \int_E \big(
		e^{i \re \biganglebracket{z,x^*}}  - 1  - i
		\kappa (z) \re \biganglebracket{z,x^*} \big)
		\lambda (dz)
		, 
	\end{align*}
        \noindent
	hence $\int_E e^{i \re \langle Tz , x^* \rangle} \mu (dz) = \int_E e^{i \re \langle z , x^* \rangle} \mu (dz)$
        for any $x^*\in E^*$, proving invariance.

        \medskip

        \noindent
   2. As the characteristic functional of $\mu$
	also 	                        takes the form \eqref{fjkldf34} 
	with
	$\displaystyle \xi (dz) = \sum_{n=0}^\infty \delta_{e_n}(dz)$
	and $\rho(e_n , du) = \delta_{\lambda_n} (du)$, $n\geq 0$,
	we note that by \eqref{e1} in Lemma~\ref{fkjdf43}
                  with $p\in [1,2)$,
	for any $c>\max_{n\geq 0} \lambda_n$
	we have
	\begin{align*}
		\bigabs{C_\mu^= (x^*, T^{*n}y^*)} & \leq
				16 \sum_{l=0}^\infty \bigbracket{2^{-p}\int_{-c}^c |u|^p |\biganglebracket{e_l,x^*}\biganglebracket{e_l,T^{*n}y^*}|^{p/2} \rho (e_l,du) + \rho (e_l,{\R\backslash[-c, c]})}
		\\
		& =
		2^{4-p}
		\sum_{l= 0}^\infty
		\lambda_l^p  |\biganglebracket{e_l,x^*}\biganglebracket{e_l,T^{*n}y^*}|^{p/2}
		\\
		&= 2^{4-p}
		\sum_{l=0}^\infty
		\lambda_{l+n}^p|\biganglebracket{e_{l+n},x^*}\biganglebracket{e_l,y^*}\omega_{l+n-1}\dots\omega_l|^{p/2}
		\\
		&\leq 2^{4-p}
		\|x^*\|^{p/2}\|y^*\|^{p/2}
		\sum_{l=0}^\infty
		\lambda_{l+n}^p\prod_{i=l}^{l+n-1} \omega_{l+n}^{p / 2}
		\\
		&= 2^{4-p}
		\|x^*\|^{p/2}\|y^*\|^{p/2}
		\sum_{l=0}^\infty 
		\lambda_{l+n}^p\prod_{j=l}^{l+n-1} \bigbracket{\frac{\lambda_j}{\lambda_{j+1}}}^{p/2}
		\\
		&= 2^{4-p}
		\|x^*\|^{p/2}\|y^*\|^{p/2}
		\sum_{l=0}^\infty 
		\lambda_{l+n}^p \bigbracket{\frac{\lambda_l}{\lambda_{l+n}}}^{p/2}
		\\
		&= 2^{4-p}
		\|x^*\|^{p/2}\|y^*\|^{p/2}
		\sum_{l=0}^\infty
		( \lambda_l \lambda_{l+n} )^{p/2}, \quad x^* \in E^*. 
	\end{align*}
        In order to conclude from Theorem~\ref{fkds3}
        it suffices to note that
        $x^* \in E^*$ is arbitrary and
        $\R\backslash Z_1 ( \nu_{x^*} )$
        is never empty by
        Lemma~\ref{genscaling-0}. 
	\end{Proof}
\noindent
 In the framework of Proposition~\ref{poiprop},
 the random variable
 $\langle X,x^*\rangle$
 has the distribution of
 $\sum_{n\geq 0} \lambda_n N_n e_n$,
 where $(N_n)_{n\geq 0}$ is a sequence of
 independent standard Poisson random variables. 
 The L\'evy measure of $\langle X,x^*\rangle$ is 
$$
 \nu_{x^*} =
 \sum_{n\geq 0} \delta_{\lambda_n \langle e_n , x^*\rangle}
 $$
 on $\real$, and we have 
$$
 Z_1 ( \nu_{x^*} )
 = \left\{ \frac{2\pi k}{\lambda_n \langle e_n , x^*\rangle}
 \ : \ 
  n\geq 0,
  \ k \in \Z
  \right\}. 
$$ 
 In particular, for $x^*\in E^*$ such that
 $\lambda_0 \langle e_0 , x^*\rangle = 2\pi$, 
 $\nu_{x^*}$ does not satisfy Condition~\ref{fjklf34-0} and 
 Proposition~2.2 of \cite{MP23}
 does not apply, as in the next example.
\begin{example}
	\label{poicor}
	Let $p\in[1,2)$, $\gamma > 1$, and let
	$T$ be the bounded weighted backward shift operator defined as 
	\begin{align*}
		Te_0 := 0, \quad Te_1 = e_0, \quad Te_{n+1} :=
		\bigg( 1 + \frac{1}{n}\bigg)^{\gamma / p}
		e_n,\qquad n\geq 1,
	\end{align*}
	and consider the compound Poisson measure $\mu$ on
        the (real) sequence space $E =\ell^p(\N)$
	with L\'evy measure \eqref{fjk2},
	where
	$\lambda_n: = \lambda_0 / (n+1)^{\gamma / p}$,
	$n\geq 0$, 
	for some $\lambda_0>0$. Then,
	$T$ admits $\mu$ as invariant measure,
	and it is mixing with respect to $\mu$,
	with the rate 
\begin{empheq}[left={\displaystyle \sup_{x^*,y^*\in E^*\setminus \{0\}}
\frac{\bigabs{C_\mu^= (x^*, T^{*n}y^*)}}{\|x^*\|^{p/2}\|y^*\|^{p/2}} \leq \empheqlbrace}]{align}
\nonumber 
& 2^{4-p} \lambda_0^p {\rm B} \left( 1 - \frac{\gamma}{2} , \gamma - 1 \right) n^{-(\gamma-1)}, \quad 1 < \gamma < 2,
\\[0.5ex]
  \label{jkldsf11} 
  & \raisebox{0.5ex}{$\displaystyle 2^{4-p} \lambda_0^p {\rm B} \left(\frac{\eps}{2} , 1 -\eps \right)n^{-(1-\eps)},
  \quad \quad \quad \gamma \geq 2,$} 
\end{empheq}
	$n\geq 1$, for any $\eps \in (0,1)$
        in \eqref{jkldsf11}, where ${\rm B}(\cdot,\cdot)$ denotes the beta function.
\end{example} 
\begin{Proof}
 For any $\gamma \in (1,2)$, we	have 
	\begin{align}
          \nonumber
	  \sum_{l=0}^\infty
	( \lambda_l \lambda_{l+n}  )^{p/2} 
	& = \sum_{l=0}^\infty \frac{\lambda_0^p}{(l+n+1)^{\gamma / 2} (l+1)^{\gamma / 2}}
	 \\
\nonumber         & \leq 
          \int_0^\infty \frac{\lambda_0^p}{(x+n)^{\gamma /2} x^{\gamma/ 2}} \dx
	  \\
          \nonumber          &= \lambda_0^p
          \int_{0}^\infty
	  \frac{n^{1-\gamma} }{(x+1)^{\gamma /2}x^{\gamma /2}} \dx
          \\
	  \label{b1}
	  &=
          2^{4-p} \lambda_0^p {\rm B} \left( 1 - \frac{\gamma}{2} , \gamma - 1 \right) n^{-(\gamma-1)}, 
	\end{align}
	and we conclude from Proposition~\ref{poiprop}.
	In the case $\gamma \geq 2$, we observe similarly 
        that for any $\eps \in (0,1)$ we have 
	\begin{align}
          \nonumber
          \sum_{l=0}^\infty
	( \lambda_l \lambda_{l+n}  )^{p/2} 
        & =	\sum_{l=0}^\infty
		\frac{\lambda_0^p}{
		  (l+n+1)^{\gamma / 2} 
	          (l+1)^{\gamma / 2}
		}
                \\
                \nonumber
		& \leq \sum_{l=0}^\infty
		\frac{\lambda_0^p}{
                  (l+n+1)^{(2 - \eps )/2}
			(l+1)^{(2 - \eps )/2}
		}
                \\
                \label{b2-2}
                & 
		\leq \lambda_0^p
                {\rm B} \left(\frac{\eps}{2} , 1-\eps \right)n^{-(1-\eps)},
	\end{align}
		and in both cases we conclude from
                                Proposition~\ref{poiprop}. 
\end{Proof} 

\section{Stable measures} 
\label{s4.2}
\noindent
In this section, we
consider the case where
$\mu$ is an 
$\alpha$-stable distribution,
$\alpha \in (0,2)\setminus \{1\}$ with characteristic functional
\eqref{st}, 
i.e. $\sigma^2 \equiv 0$ and 
$\rho (z , du)$ in \eqref{fjkldf34} 
takes the form 
$$
\rho(z,du) = \frac{du}{|u|^{1+\alpha}},\qquad u\neq 0, 
$$
see the discussion following \cite[Theorem 4]{SCMIDP}. 
In particular, in Proposition~\ref{extexp}
we will derive mixing rates for
a family of weighted shifts
that leave $\alpha$-stable measures invariant 
on the sequence space $\ell^p(\N)$,
$p\in[1,2]$. 
In Corollary~\ref{infseriesrate}
we also derive decay rates for 
quantities of the form \eqref{fjkl1} 
for finite and infinite
linear combinations 
of exponentials. 
\begin{lemma}
	\label{fkjldsf-1}
	Let $\mu$
	be an $\alpha$-stable distribution 
	with control measure $\xi$ on $E$,
	$\alpha\in (0,2) \setminus \{1\}$. 
	For any $p\in( \alpha ,2]$ and $c>0$, we have the codifference bounds 
	\begin{align*}
		\bigabs{C_\mu^= (x^*, y^*)} &\leq \frac{2^{5-p}c^{p-\alpha}}{p-\alpha}
		\int_E |\re\biganglebracket{z,x^*}\re\biganglebracket{z,y^*}|^{p/2} \xi(dz)
		+ \frac{32}{\alpha c^\alpha}\xi(E) 
	\end{align*}
	and 
	\begin{align*}
		\bigabs{C_\mu^{\neq}(x^*, y^*)} &\leq \frac{2^{5-p}c^{p-\alpha}}{p-\alpha}
		\int_E |\re\biganglebracket{z,x^*}\imag\biganglebracket{z,y^*}|^{p/2} \xi(dz)
		+ \frac{32}{\alpha c^\alpha }\xi(E) ,
	\end{align*}
	for       $x^*,y^*\in E^*$.
\end{lemma} 
\begin{Proof}
	  By Lemma~\ref{fkjdf43}, for
  $p\in (\alpha,2]$ we have
	\begin{align*}
		\bigabs{C_\mu^= (x^*, y^*)}
			& \leq 16
		\hskip-0.06cm
		\int_E
		\hskip-0.1cm
		\bigbracket{
			2^{-p}
			|\re\biganglebracket{z,x^*}\re\biganglebracket{z,y^*}|^{p/2} \int_{-c}^c \frac{1}{|u|^{1+\alpha-p}} \du + \int_{\R\backslash[-c, c]} \frac{1}{|u|^{1+\alpha}} \du}\xi(dz)
                \\
                & =
		16 \int_E
		\hskip-0.1cm
		\bigbracket{
		  \frac{
                    2^{1-p}
                    c^{p-\alpha}}{p-\alpha}
			|\re\biganglebracket{z,x^*}\re\biganglebracket{z,y^*}|^{p/2}
 + \frac{2}{\alpha c^\alpha} }\xi(dz)
 \\
 & =
 \frac{2^{5-p}c^{p-\alpha}}{p-\alpha}
		\int_E |\re\biganglebracket{z,x^*}\re\biganglebracket{z,y^*}|^{p/2} \xi(dz)
		+ \frac{32}{\alpha c^\alpha }\xi(E). 
        \end{align*}
        The proof is similar for $C_\mu^{\neq}(x^*, y^*)$.
\end{Proof}
In the next proposition,
for $p\in (\alpha , 2]\cap [1,2]$ (so that $\ell^p(\Z)$ is a Banach space and Lemma~\ref{fkjldsf-1} applies),
we provide 
rates for the 
mixing of weighted forward shift operators on
$\ell^p(\Z)$
considered in Proposition~4.2 of \cite{MP23}.  
In what follows, we write $f(n) = O_y(g(n))$ when we have
$$
| f(n) | \leq C_y |g(n)|, \quad n \geq 0, 
$$
for $C_y>0$ a constant possibly depending on $y$ and independent of $n\geq 0$.
Also, as above, we let $(e_n)_{n\geq 0}$ denote the canonical basis
of $\ell^p ( \N )$.
\begin{prop}
\label{extexp}
	Let  $\alpha \in (0,2)\setminus \{1\}$ 
	and $p\geq 1$. Let $(\omega_n)_{n\in \Z}$ be a
	positive weight sequence
	such that the sequence $(k_n)_{n\in\Z}$ defined by
	$$
	k_n := 
	k_0
	{\bf 1}_{\{n\leq -1\}} 
	\prod_{l=n+1}^0 \frac{1}{\omega_l}
	+
	k_0
	{\bf 1}_{\{n\geq 0\}} 
	\prod_{l=1}^n \omega_l.
	$$
        belongs to $\ell^\alpha(\Z)$. 
	Let $\mu$ be the $\alpha$-stable measure
	on $E=\ell^p (\Z )$, $p\geq 1$, 
	with characteristic functional \eqref{st}
	and control measure
	\begin{align}
			\nonumber 
			\xi (dz) := \frac{1}{2} 
			\sum_{n=-\infty}^\infty 
			k_n^\alpha
			\big(
			\delta_{e_n} (dz) + \delta_{ie_n} (dz)
			\big)
			,
	\end{align}
	and 
	consider the weighted forward shift operator $T$ 
	on $E$
	defined by
	$$Te_n := \omega_{n+1}e_{n+1}, \quad n\in \Z.
	$$
	The following are true.
	\begin{enumerate}
		\item $T$ admits $\mu$ as invariant measure.
		\item Assume that $p\in (\alpha , 2]\cap [1,2]$
		and              
		there exist 
		$q_-\geq 0$, $q_+\geq 1$ such that
		$$
		\eta_- : = \sup_{l\leq -q_-} \frac{1}{\omega_l} < 1 
		\quad
		and
		\quad
		\eta_+ : = \sup_{l\geq q_+} \omega_l < 1
		$$
		with $\eta_+^{p/2} \neq \eta_-^{\alpha - p/2}$.
				Then,
				\begin{equation}  \label{basicest}
			\sup_{x^*,y^*\in E^*\setminus \{0\}}
			\frac{
				\bigabs{C_\mu^{=, \neq}
					(x^*, T^{*n}y^*)}
			} {                  	        \|x^*\|^{\alpha /2} \|y^*\|^{\alpha /2}}
			=
			O_{\mu ,T}\big(
			\max \big(
			\eta_-^{2\alpha/p - 1}
			,
			\eta_+
			\big)^{\alpha n/2}\big),
			\quad
			n \geq 0. 
		\end{equation}
		In particular, $T$
		is mixing when $\alpha \in (1/2 , 2) \setminus \{1 \}$ and $p\in (\alpha , 2\alpha) \cap [1,2]$.
	\end{enumerate}
\end{prop}
\begin{Proof} 
\noindent
  1. The condition $\sum_{n=-\infty}^\infty |k_n|^\alpha < \infty$ ensures that $\xi$ is finite as a control measure.
	To show that $T$ admits $\mu$ as invariant measure, we
	use \eqref{st} and the equalities 
	\begin{align*} 
		\int_E \bigabs{\re  \biganglebracket{Tz,x^*}}^\alpha \xi (dz)
		& = 
		\sum_{n=-\infty}^\infty
		k_n^\alpha \bigabs{\re  \biganglebracket{Te_n,x^*}}^\alpha
		\\
		& = 
		\sum_{n=-\infty}^\infty  
		k_n^\alpha\omega_{n+1}^\alpha \bigabs{\re  \biganglebracket{e_{n+1},x^*}}^\alpha
		\\
		& = 
		\sum_{n=-\infty}^\infty  
		k_{n+1}^\alpha \bigabs{\re  \biganglebracket{e_{n+1},x^*}}^\alpha
				\\
		& =
		\int_E \bigabs{\re  \biganglebracket{z,x^*}}^\alpha \xi (dz),
	\end{align*} 
	hence $\int_E e^{i \re \langle Tz , x^* \rangle} \mu (dz) = \int_E e^{i \re \langle z , x^* \rangle} \mu (dz)$
        for all $x^*\in E^*$, proving invariance.

        \medskip
        
        \noindent
  2. Without loss of generality, we may assume
	$q_-=0$ and $q_+=1$.
	We note that the right-hand side term in
	Lemma~\ref{fkjldsf-1} can be bounded as 
	\begin{align}
		\nonumber
		\int_E |\re\biganglebracket{z,x^*}\re\biganglebracket{z,T^{*n}y^*}|^{p/2} \xi(dz)
		& \leq
		\|x^*\|^{p/2}
		\int_{E}
		\|z\|^{p/2}
		|\re\biganglebracket{z,T^{*n}y^*}|^{p/2}
		\xi(dz)
		\\
		\nonumber
		&  = 
		\|x^*\|^{p/2}
		\int_{E}
		\|z\|^{p/2}|\re\biganglebracket{z,T^{*n}y^*}|^{p/2}
		\xi(dz)
		\\
		\nonumber
		&  = 
		\|x^*\|^{p/2}
		\int_{E}
		\|z\|^{p/2}|\re\biganglebracket{T^nz,y^*}|^{p/2}
		\xi(dz)
		\\
		\nonumber 
		&  \leq 
		\|x^*\|^{p/2}
		\|y^*\|^{p/2}
		\int_{E}
		\|z\|^{p/2}\|T^nz\|^{p/2}
		\xi(dz)
		, 
	\end{align}
	$x^*,y^*\in E^*$, $n \geq 0$.
	We have
	\begin{align}
		\nonumber
		\int_{E}
		\|z\|^{p/2}\|T^nz\|^{p/2}
		\xi(dz) &= 
		\frac{1 
		}{2} 
		\sum_{l=-\infty}^\infty 
		k_l^\alpha
		\int_E
		\|z\|^{p/2}\|T^nz\|^{p/2}
		\big(
		\delta_{e_l} (dz) + \delta_{ie_l} (dz)
		\big)
		\\
		&
		\nonumber  
		= 
				\sum_{l=-\infty}^\infty k_l^\alpha
		\|T^n e_l\|^{p/2}
		\\
		\nonumber 
		&=
				\sum_{l=-\infty}^\infty k_l^\alpha \prod_{j=l+1}^{l+n} \omega_j^{p/2}.
	\end{align}
		We split the above series into three components. 
	\begin{itemize}
		\item If $l\geq 0$, then
		$$
		\prod_{j=l+1}^{l+n} \omega_j^{p/2}\leq \eta_+^{pn/2},
		$$
		and so
		\begin{align*}
			\sum_{l=0}^\infty k_l^\alpha \prod_{j=l+1}^{l+n} \omega_j^{p/2}
			\leq k_0^\alpha \sum_{l=0}^\infty \eta_+^{\alpha l} \eta_+^{pn/2}
			=   \frac{k_0^\alpha}{1-\eta_+^\alpha}\eta_+^{pn/2}.
		\end{align*}
		\item If $-n < l \leq -1$, then
		$$
		k_l^\alpha \prod_{j=l+1}^{l+n} \omega_j^{p/2} =
		k_0^\alpha
		\prod_{j=l+1}^0 \frac{1}{\omega_j^{\alpha - p/2}} \prod_{j=1}^{l+n} \omega_j^{p/2} \leq
		k_0^\alpha \eta_-^{(\alpha - p/2)|l|} \eta_+^{(l+n)p/2},
		$$
		and so
		\begin{align*}
			\sum_{l=-n+1}^{-1} k_l^\alpha \prod_{j=l+1}^{l+n} \omega_j^{p/2} &\leq
			k_0^\alpha
			\sum_{l=-n+1}^{-1} \eta_-^{-(\alpha - p/2)|l|} \eta_+^{(l+n)p/2}\\
			&=
			k_0^\alpha
			\eta_+^{pn/2}\sum_{l=-n+1}^{-1} \eta_-^{-(\alpha - p/2)l} \eta_+^{pl/2}\\
			&\leq
			k_0^\alpha
			\eta_+^{pn/2} \sum_{l=1}^{n+1} \bigbracket{\frac{\eta_-^{\alpha - p/2}}{\eta_+^{p/2}}}^l\\
			&=
			k_0^\alpha
			\eta_+^{pn/2}
                        \bigbracket{\frac{\eta_-^{\alpha - p/2}}{\eta_+^{p/2}}} \frac{1 - \bigbracket{\frac{\eta_-^{\alpha - p/2}}{\eta_+^{p/2}}}^{n+1}}{1-\bigbracket{\frac{\eta_-^{\alpha - p/2}}{\eta_+^{p/2}}}}\\
			&= \frac{k_0^\alpha\eta_-^{\alpha - p/2}}{\eta_+^{p/2}-\eta_-^{\alpha - p/2}}\eta_+^{pn/2}
			-
			\frac{k_0^\alpha\eta_-^{\alpha - p/2}}{\eta_+^{p/2}-\eta_-^{\alpha - p/2}}\bigbracket{\frac{\eta_-^{\alpha - p/2}}{\eta_+^{p/2}}}\eta_-^{(\alpha - p/2)n}.
		\end{align*}
		\item If $l\leq -n$, then
		$$
		k_l^\alpha \prod_{j=l+1}^{l+n} \omega_j^{p/2}
		= k_0^\alpha
		\prod_{j=l+1}^0 \frac{1}{\omega_j^{\alpha - p/2}} \prod_{j=l+n+1}^{0} \frac{1}{\omega_j^{p/2}} \leq k_0^\alpha \eta_-^{(\alpha - p/2)|l|} \eta_-^{p|l+n|/2},
		$$
		and so 
		\begin{align*}
			\sum_{l=-\infty}^{n} k_l^\alpha \prod_{j=l+1}^{l+n} \omega_j^{p/2} \leq
			k_0^\alpha
			\sum_{l=n}^\infty \eta_-^{(\alpha - p/2)\l} \eta_-^{p(l+n)/2}
			= k_0^\alpha 
			 \eta_-^{pn/2}\sum_{l=n}^\infty \eta_-^{\alpha l} 
			 = \frac{k_0^\alpha}{1-\eta_-^\alpha }\eta_-^{(\alpha + p/2)n}.
		\end{align*}
	\end{itemize}
	Hence, 
	we have 
	\begin{equation}
		\nonumber 
		\int_E |\re\biganglebracket{z,x^*}\re\biganglebracket{z,T^{*n}y^*}|^{p/2} \xi(dz)
				\leq 
		k_0^\alpha
		\|x^*\|^{p/2} \|y^*\|^{p/2} 
		\big(
		K_1 \eta_+^{pn/2} + K_2\eta_-^{(\alpha + p/2)n} + K_3\eta_-^{(\alpha - p/2)n}
		\big), 
	\end{equation}
	$x^*,y^*\in E^*$,
	where
	$$
	K_1:=                                 
	\frac{1}{1-\eta_+^\alpha} + \frac{\eta_-^{\alpha - p/2}}{\eta_+^{p/2}-
          \eta_-^{\alpha - p/2}},
	\quad 
	K_2:= 
	\frac{1}{1-\eta_-^\alpha },
	\quad 
	K_3:= 
	-\frac{\eta_-^{\alpha - p/2}}{\eta_+^{p/2}-\eta_-^{\alpha - p/2}}
	\bigbracket{\frac{\eta_-^{\alpha - p/2}}{\eta_+^{p/2}}}
	$$
	are constants independent of $n \geq 1$.
		Observe that since $\eta_- < 1$, we have $\eta_-^{(\alpha + p/2)n} = O(\eta_-^{(\alpha - p/2)n})$. If $\eta_+^{p/2} > \eta_-^{\alpha - p/2}$ then $K_3 < 0$ and
				\begin{equation}
				\label{fjkdl12}
				\int_E |\re\biganglebracket{z,x^*}\re\biganglebracket{z,T^{*n}y^*}|^{p/2} \xi(dz)
				=  \|x^*\|^{p/2} \|y^*\|^{p/2}
				O_{\mu ,T} \big(\eta_+^{pn/2}\big).
			\end{equation}
			Conversely, if $\eta_+^{p/2} < \eta_-^{\alpha - p/2}$ then $K_3 > 0$. We have $\eta_+^{pn/2} = O(\eta_-^{(\alpha - p/2)n})$, and thus 
			\begin{equation} 
				\label{fjkdl13}
				\int_E |\re\biganglebracket{z,x^*}\re\biganglebracket{z,T^{*n}y^*}|^{p/2} \xi(dz)
				=
				\|x^*\|^{p/2} \|y^*\|^{p/2}
				O_{\mu ,T} \big(\eta_-^{(\alpha - p/2)n}\big).
			\end{equation}
		We now bound the codifferences. 
	If $\eta_+^{p/2} > \eta_-^{\alpha - p/2}$, since $p\in( \alpha ,2]$, 	by Lemma~\ref{fkjldsf-1} and \eqref{fjkdl12} 
	for any $c>0$ we have 
	\begin{align*}
		\bigabs{C_\mu^= (x^*, T^{*n}y^*)}
		=  \|x^*\|^{p/2} \|y^*\|^{p/2}
		O_{\mu ,T} \big(c^{p-\alpha}\eta_+^{pn/2}\big) + O_\mu (c^{-\alpha}).
	\end{align*}
	In particular, letting $c_n :=
	\|x^*\|^{-1/2} \|y^*\|^{-1/2} \eta_+^{-n/2}$,
	$n \geq 1$, putting $c = c_n$ we get
	\begin{align*}
		\bigabs{C_\mu^= (x^*, T^{*n}y^*)}
		= \|x^*\|^{\alpha /2} \|y^*\|^{\alpha /2} O_{\mu ,T} \big(\eta_+^{\alpha n/2}\big). 
	\end{align*}
	In the case
	$\eta_+^{p/2} < \eta_-^{\alpha - p/2}$, a similar argument
	using \eqref{fjkdl13}
	with
	$c_n :=
	\|x^*\|^{-1/2} \|y^*\|^{-1/2} \eta_-^{(1/2 - \alpha/p)n}$,
	$n \geq 1$,
	gives
	\begin{align*}
		\bigabs{C_\mu^= (x^*, T^{*n}y^*)}
		= \|x^*\|^{\alpha /2} \|y^*\|^{\alpha /2} O_{\mu ,T} \big(
                \eta_-^{\alpha(\alpha/p - 1/2)n}\big). 
	\end{align*}
	Likewise, a similar argument establishes the corresponding inequality for $\bigabs{C_\mu^{\neq} (x^*, T^{*n}y^*)}$,
	hence we have
	\begin{align}
          \nonumber
	  \bigabs{C_\mu^{=, \neq}
			(x^*, T^{*n}y^*)}
		=
		\left\{ 
		\begin{array}{ll}
			\|x^*\|^{\alpha /2} \|y^*\|^{\alpha /2}
			O_{\mu ,T}\big(\eta_+^{\alpha n/2}\big)
			& \mbox{if }
			\eta_+^{p/2} > \eta_-^{\alpha - p/2},
			\bigskip
			\\  
			\|x^*\|^{\alpha /2} \|y^*\|^{\alpha /2}
			O_{\mu ,T}\big(
                        \eta_-^{
				\alpha(\alpha/p - 1/2)n
			}\big)
			& \mbox{if }  
			\eta_+^{p/2} < \eta_-^{\alpha - p/2}, 
		\end{array}
		\right.
	\end{align}
  which rewrites as
  \begin{align*}
	\bigabs{C_\mu^{=, \neq}
		(x^*, T^{*n}y^*)}
	=\|x^*\|^{\alpha /2} \|y^*\|^{\alpha /2}
	O_{\mu ,T}\big(\max\big(\eta_-^{2\alpha/p - 1}, \eta_+\big)^{\alpha n/2} \big),
\end{align*}
 giving \eqref{basicest}. Mixing in the case of $\alpha > 1/2$
 and $p\in (\alpha , 2\alpha) \cap [1,2]$ follows from
 Theorem~\ref{fkds3} and the decay of codifferences.
\end{Proof}
 Proposition~\ref{extexp} also applies to
  symmetric control measures of the form 
		\begin{align}
			\nonumber 
			\xi (dz) := \frac{1}{4} 
			\sum_{n=-\infty}^\infty 
			k_n^\alpha
			\big(
			\delta_{e_n} (dz) + \delta_{-e_n} (dz) + \delta_{ie_n} (dz) + \delta_{-ie_n} (dz)
			\big). 
		\end{align} 
\noindent
We now extend Proposition~\ref{extexp}
by determining decay rates for the quantity 
\begin{align*}
I_n(f,g) := \int_E f(z)g(T^n z) \mu(dz) - \int_E f(z)\mu(dz) \int_E g(z)\mu(dz), \quad n\geq 0, 
\end{align*}
for a family of functions $f,g$ in $L^2(E,\mu)$, when $T$ is mixing. 
Table~\ref{fctable} displays the equivalent codifferences
\begin{equation}
\nonumber 
C^{\phi , \psi}_\mu (x^*, y^*) := \log
\int_E
e^{i \phi (\langle z,x^*\rangle) + i\psi (\langle z,y^*\rangle)} 
\mu (dz)
- \log
\int_E
e^{i \phi (\langle z,x^*\rangle)}\mu (dz)
- \log
\int_E
e^{i \psi (\langle z,y^*\rangle)}
\mu (dz)
,
\end{equation} 
$x^*,y^*\in E^*$,
for different choices of functions $\phi$, $\psi$. 
\begin{table}[H]
\centering
\begin{tabular}{|c|c||c|}
	\hline
	$\phi(\cdot)$ & $\psi(\cdot)$ & $C^{\phi , \psi}_\mu (x^*, y^*)$ \\ 
	\hline
	$\re(\cdot)$ & $-\re(\cdot)$ & $C_\mu^= (x^*, y^*)$\\
	\hline
	$-\re(\cdot)$ & $\re(\cdot)$ & $C_\mu^= (-x^*, -y^*)$\\
	\hline
	$\imag(\cdot)$ & $-\imag(\cdot)$ & $C_\mu^= (-ix^*, -iy^*)$\\
	\hline
	$-\imag(\cdot)$ & $\imag(\cdot)$ & $C_\mu^= (ix^*, iy^*)$\\
	\hline
	$\re(\cdot)$ & $-\imag(\cdot)$ & $C_\mu^{\neq}(x^*, y^*)$\\
	\hline
	$-\re(\cdot)$ & $\imag(\cdot)$ & $C_\mu^{\neq}(-x^*, -y^*)$\\
	\hline
	$\imag(\cdot)$ & $\re(\cdot)$ & $C_\mu^{\neq}(-ix^*, -iy^*)$\\
	\hline
	$-\imag(\cdot)$ & $-\re(\cdot)$ & $C_\mu^{\neq}(ix^*, iy^*)$\\
	\hline
	$-\imag(\cdot)$ & $\re(\cdot)$ & $C_\mu^{\neq}(x^*, y^*)$\\
	\hline
	$\imag(\cdot)$ & $-\re(\cdot)$ & $C_\mu^{\neq}(-x^*, -y^*)$\\
	\hline
	$\re(\cdot)$ & $\imag(\cdot)$ & $C_\mu^{\neq}(-ix^*, -iy^*)$\\
	\hline
	$-\re(\cdot)$ & $-\imag(\cdot)$ & $C_\mu^{\neq}(ix^*, iy^*)$\\
	\hline
\end{tabular}
\caption{Function-codifference triples}
\label{fctable}
\end{table}

\vspace{-0.3cm}

\noindent
First, we observe that the estimate \eqref{basicest} 
holds for the twelve codifference quantities in Table~\ref{fctable}.
 Next, we derive bounds on
$I_n(f,g)$ for
$f$ and $g$ in a class
of functions defined by 
infinite series. 
\begin{corollary} 
\label{infseriesrate}
Let $\mu$ be the $\alpha$-stable measure
on $E=\ell^p (\Z )$
for 
$\alpha\in (1/2,2) \setminus \{1\}$ and $p\in (\alpha , 2\alpha) \cap [1,2]$,        let $T$
be the mixing weighted forward shift operator 
defined in Proposition~\ref{extexp},
 let $(\phi , \psi)$ be 
a pair of functions given in Table~\ref{fctable}, and 
let $(a_j)_{j\in \N}$ and $(b_l)_{l\in \N}$
be complex $\ell^1(\N)$ sequences such that
$$
\sum_{j=0}^\infty |a_j|\|T\|^{jp/2} < \infty
\quad
\mbox{and}
\quad
\sum_{l=0}^\infty |b_l|\|T\|^{lp/2} < \infty.
$$ 	
 For any $x^*,y^*$, the functions 
\begin{align*}
	f (z) := \sum_{j=0}^\infty a_j e^{i\phi
		( \langle z,T^{*j}x^*\rangle )}
	\ \ \text{and} \ \
	g (z) := \sum_{l=0}^\infty b_l e^{i\psi
		( \langle z,T^{*l}y^* \rangle )}
\end{align*}
are well-defined in $L^2(E,\mu)$, and we have 
\begin{equation}
	\label{rate1} 
	\bigabs{I_n\bigbracket{f,g}}
	=
	O_{f,g,\mu ,T}\big(\max \big(
	\eta_-^{2\alpha/p - 1}
	,
	\eta_+
	\big)^{\alpha n/2}\big), \qquad n \geq 0. 
\end{equation} 
\end{corollary} 
\begin{Proof}
For any pair $(\phi , \psi )$  of functions in
Table~\ref{fctable}, we have 
\begin{align*}
	\big|
	I_n\big(e^{i\re \phi ( \biganglebracket{\cdot,x^*} )},
	e^{-i\re\psi ( \biganglebracket{\cdot,y^*})}\big) 
	\big|
	&
	=
	\big|
	\exp(C^{\phi,\psi}_\mu (x^*, y^*)) - 1
	\big|
	\left|
	\int_E e^{i\re\biganglebracket{z,x^*}} \mu(dz)
	\int_E e^{-i\re\biganglebracket{z,y^*}} \mu(dz)
	\right|
	\\
	&\leq 
	\big| \exp(C^{\phi,\psi}_\mu (x^*, y^*)) - 1\big|
	\\
	&       \leq 
	|C_\mu^{\phi , \psi } ( x^*, y^*)|
	\exp\left( \big|C_\mu^{\phi , \psi }\big( x^*, y^*\big)\big|\right). 
\end{align*} 
We have, if $0 \leq n \leq j-l$, 
$$             \exp\left( \big|C_\mu^{\phi , \psi } \big(T^{*j}x^*, T^{*(n+l)}y^*\big)\big|\right) =
\exp\left(
\big|C_\mu^{\phi , \psi } \big(T^{*(j-n-l)}x^*, y^*\big)\big|
\right) 
\leq 
\exp\left(
\max_{n\geq 0} \big|C_\mu^{\phi , \psi } \big(T^{*n}x^*, y^*\big)\big|
\right), 
$$
and, if $n\geq  \max ( 0,j-l)$, 
$$
\exp\left( \big|C_\mu^{\phi , \psi } \big(T^{*j}x^*, T^{*(n+l)}y^*\big)\big|\right) 
=  
\exp\left( \big|C_\mu^{\phi , \psi } \big(x^*, T^{*(n+l-j)}y^*\big)\big|\right)
\leq 
\exp\left(
\max_{n\geq 0} \big|C_\mu^{\phi , \psi } \big(x^*, T^{*n}y^*\big)\big|
\right), 
$$
where the maxima are finite by
\eqref{basicest}.                 
Hence, 
letting 
$$
f_{j,x^*}(z) := e^{i\phi
	(\langle z,T^{*j}x^* \rangle )}
\quad
\mbox{and}
\quad
g_{l,y^*}(z) := e^{i\psi
	(\langle z,T^{*l}y^*\rangle )},
\quad j,l \geq 0, 
$$
 by Proposition~\ref{extexp} we have
\begin{align*} 
	&\bigabs{I_n\bigbracket{f_{j,x^*},g_{l,y^*}}}\\
	&\quad\leq 
	|C_\mu^{\phi , \psi } (T^{*j}x^*, T^{*(n+l)}y^*)|
	\exp\left(
	\max\left(
	\max_{n\geq 0} \big|C_\mu^{\phi , \psi } \big(T^{*n}x^*, y^*\big)\big|
	,
	\max_{n\geq 0} \big|C_\mu^{\phi , \psi } \big(x^*, T^{*n}y^*\big)\big|
	\right) \right) 
	\\
	&
	\quad             
	\leq 
	\|T^{*j}x^*\|^{\alpha /2} \|T^{*l}y^*\|^{\alpha /2}
	O_{x^*,y^*,\mu ,T} \big(\max \big(
	\eta_-^{2\alpha/p - 1}
	,
	\eta_+
	\big)^{\alpha n/2}
	\big), \quad
      j,l \geq 0. 
\end{align*} 
 Since
$(a_j)_{j\geq 0} \in \ell^1(\N)$
and 
$(b_l)_{l\geq 0}\in \ell^1(\N)$,
we have $f,g\in L^2(E,\mu)$, 
the series 
$\sum_{j=0}^\infty f_{j,x^*}(z)$ and
$\sum_{l=0}^\infty g_{l,y^*}(z)$
converge absolutely for every $z \in E$, and
\begin{align*} 
	|I_n(f,g)| & \leq
	\sum_{j,l=0}^{\infty}
	|a_j||b_l| |I_n(f_{j,x^*}, g_{l,y^*})|
	\\
	&
	\leq
	\|x^*\|^{\alpha /2}
	\|y^*\|^{\alpha /2}
	\sum_{j,l=0}^{\infty}|a_j||b_l| \|T^j\|^{\alpha /2}
	\|T^l\|^{\alpha /2} 
	O_{x^*,y^*,\mu ,T} \big(\max \big(
	\eta_-^{2\alpha/p - 1}
	,
	\eta_+
	\big)^{\alpha n/2}
	\big)
	\\
	&
	\leq
	\|x^*\|^{\alpha /2} \|y^*\|^{\alpha /2}
	O_{x^*,y^*,\mu ,T}\big(\max \big(
	\eta_-^{2\alpha/p - 1}
	,
	\eta_+
	\big)^{\alpha n/2}
	\big)                  \sum_{j,l=0}^{\infty}|a_j||b_l|
	\|T^j\|^{\alpha /2}
	\|T^l\|^{\alpha /2} 
	.
\end{align*} 
\end{Proof}
\noindent
Since Table~\ref{fctable} contains every
possible pair $(\phi,\psi) \in \{\pm\re(\cdot), \pm\imag(\cdot)\}^2$ where $\phi\neq\psi$,
the rate obtained in \eqref{rate1}
also applies
to
$f$ and $g$ given by
\begin{align*}
f(z) := \sum_{\phi\in\Phi}\sum_{j=0}^\infty a_{\phi , j} e^{i\phi(
	\langle z,T^{*j}x^* \rangle)}
\ \  \text{and} \ \ 
g(z) := \sum_{\psi\in\Psi}
\sum_{l=0}^\infty b_{\psi , l}
e^{i\psi (\langle z,T^{*l}y^* \rangle)}, 
\end{align*}
where $\Phi,\Psi$
are disjoint non-empty subsets of
$\{\pm\re(\cdot), \pm\imag(\cdot)\}$, and 
$(a_{\phi,j})_{j\in \N , \phi\in\Phi }$, 
$(b_{\psi,l})_{l\in \N , \psi\in\Psi}$
are complex $\ell^1(\N)$ sequences such that
\begin{align*}
\sum_{j=0}^\infty |a_{\phi,j}|\|T\|^{jp/2} < \infty \quad
\mbox{and}
\quad
\sum_{l=0}^\infty |b_{\psi,l}|\|T\|^{lp/2} < \infty, \quad
\phi\in\Phi, \psi\in\Psi. 
\end{align*}
\section{Tempered stable measures} 
\label{s6}
\noindent 
  In this section, we consider the case
where $\mu$ is the distribution of an $\ell^p(\N)$-valued random variable of the form
\begin{align}\label{cpvr}
	\sum_{n=0}^\infty
	k_n\big(
	\theta_{1,n} + i\theta_{2,n}
	\big)e_n,
\end{align}
where
$k_n$ is an appropriately chosen positive sequence, and $\theta_{1,n}$ and $\theta_{2,n}$ are independent and identically distributed copies of a real-valued tempered stable random variable $\theta$, centered at zero, with two-sided index of stability $\alpha \in (0,1)$ and characteristic function 
\begin{align*} 
	\E[e^{it\theta}] =  \exp \left(\int_{\R} \big(
	e^{i tx}  - 1 - i
	\kappa (x) tx
	\big)
	\lambda (dx) \right),\qquad t\in\R,
\end{align*}
where
\begin{align*}
	\kappa (x) = {\bf 1}_{\{ |x| < 1 \}}
	+ \frac{1}{|x|}
	{\bf 1}_{\{|x|\geq 1\}}, \qquad
	x \in \R,
\end{align*}
and
\begin{align}\label{lefy}
	\lambda(dx) = a_-\frac{e^{-\lambda_- |x|}}{|x|^{1+\alpha}} \mathbf{1}_{\R^-}(x) \dx
	+
	a_+\frac{e^{-\lambda_+ x}}{x^{1+\alpha}} \mathbf{1}_{\R^+}(x) \dx,
\end{align}
with $a_-, a_+, \lambda_-, \lambda_+ > 0$, see for example
\cite{koponen, TSDP}. We first derive criteria under which
the series \eqref{cpvr} is well-defined, i.e. is almost surely $\ell^p(\N)$-valued, as in \cite{Owls, MP23}. 
\begin{prop}
	\label{tsconv}
	Let $\alpha\in (0,1)$ and 
	$(k_n)_{n\geq 0}\in\ell^\alpha (\N)$. 
         The following are true.
	\begin{enumerate}
		\item The series $\disp \sum_{n=0}^\infty |k_n\theta_{1,n}|^p$
		converges almost surely 
		for all $p\in (\alpha,\infty)$.
	      \item The series $\disp \sum_{n=0}^\infty |k_n\theta_{1,n} + ik_n\theta_{2,n}|^p$
                 converges almost surely 
		for all $p\in [1,\infty)$.
	\end{enumerate}
\end{prop}
\begin{Proof}
  To prove the first statement, we use the three-series theorem of Kolmogorov,
  see e.g. \cite{durrett2005}.
   For notational simplicity write $\theta_n = \theta_{1,n}$. Let
	\begin{align*}
		\theta_n' = \begin{cases}
			\theta_n, & \displaystyle
			|\theta_n| < \frac{1}{|k_n|},
			\medskip
			\\
			0, & \displaystyle |\theta_n| \geq \frac{1}{|k_n|},
		\end{cases}
	\end{align*}
	and let $f(x)$ denote the common probability density of $\theta_n$,
        $n\geq 0$. 
	\begin{itemize}
		\item The first series condition requires to show that 
		\begin{align*}
			\sum_{n=0}^\infty \P(|k_n\theta_n| > 1) &= \sum_{n=0}^\infty\bigbracket{\int_{1/|k_n|}^\infty f(x)\dx + \int_{-\infty}^{-1/|k_n|} f(x)\dx}\\
			& = \sum_{n=0}^\infty\int_{1/|k_n|}^\infty f(x)\dx + \sum_{n=0}^\infty\int_{-\infty}^{-1/|k_n|} f(x)\dx\\
			&< \infty. 
		\end{align*}
                                 We show finiteness of the first series; the second is similar. From \cite[Theorem 7.10]{TSDP} we have
		$$f(x)\sim
		C x^{-1-\alpha} e^{-\lambda_+ x}
		$$
		as $x\to\infty$, where
		\begin{align*}
			C = a_+\exp\bigbracket{-a_+\Gamma(-\alpha)(\lambda_+)^\alpha + a_-\Gamma(-\alpha)\bigsquarebracket{(\lambda_+ + \lambda_-)^\alpha - (\lambda_-)^\alpha }}
		\end{align*}
		is a positive constant.
				We conclude from 
		\begin{align*}
			\int_{1/|k_n|}^\infty f(x)\dx \sim C \int_{1/|k_n|}^\infty \frac{e^{-\lambda_+ x}}{x^{1+\alpha}} \dx \leq C \int_{1/|k_n|}^\infty \frac{1}{x^{1+\alpha}} \dx = \frac{C}{\alpha}|k_n|^\alpha
		\end{align*}
	and the fact that 
			$(k_n)_{n\geq 0}\in\ell^\alpha (\N)$.
        			\item The second series condition requires to show that 
		\begin{align*}
			\sum_{n=0}^\infty \E[|k_n\theta_n'|^p] &= \sum_{n=0}^\infty |k_n|^p \bigbracket{\int_{1/|k_n|}^\infty x^pf(x)\dx + \int_{-\infty}^{-1/|k_n|} x^pf(x)\dx}\\
			&= \sum_{n=0}^\infty |k_n|^p \int_{1/|k_n|}^\infty x^pf(x)\dx + \sum_{n=0}^\infty |k_n|^p\int_{-\infty}^{-1/|k_n|} x^pf(x)\dx\\
			&< \infty.
		\end{align*}
		Using the asymptotics 
                $
			f(x) \sim
			C x^{-1-\alpha} e^{-\lambda_+ x}
		$
		as $x$ tends to infinity, 
		since $p > \alpha$ we have 
		\begin{align*}
 \int_1^{k/|k_n|} x^p f(x)\dx &\sim C \int_0^{1/|k_n|} \frac{e^{-\lambda_+ x}}{x^{1+\alpha-p}} \dx\\
			&\leq C \int_0^{1/|k_n|} \frac{1}{x^{1+\alpha-p}} \dx\\
			&= \frac{C}{\alpha - p} |k_n|^{\alpha-p},
		\end{align*}
		hence 	since $(k_n)_{n\geq 0}\in\ell^\alpha(\N)$,
                we have 
		\begin{align*}
			\sum_{n=0}^\infty \E[|k_n\theta_n'|^p] \leq
			\sum_{n=0}^\infty \bigbracket{|k_n|^p \cdot \frac{C}{\alpha - p} |k_n|^{\alpha-p}}
			< \infty.
		\end{align*}
			\item The third series condition  requires to show that 
		\begin{align*}
			\sum_{n=0}^\infty \Var[|k_n\theta_n'|^p] &= \sum_{n=0}^\infty \bigbracket{\E[|k_n\theta_n'|^{2p}] - \E[|k_n\theta_n'|^p]^2}\\
			&=  \sum_{n=0}^\infty\E[|k_n\theta_n'|^{2p}] - \sum_{n=0}^\infty\E[|k_n\theta_n'|^p]^2 \\
			&< \infty. 
		\end{align*}
				Observe that
		\begin{align*}
			\sum_{n=0}^\infty \E[|k_n\theta_n'|^{2p}] < \infty
		\end{align*}
		by applying the second series condition argument with $2p > \alpha$. Since variance is non-negative, and the second series has only non-negative terms, it follows that this condition also holds.
	\end{itemize}
	Finally, the second statement follows from the first by the inequality
	\begin{align*}
		|a+ib|^p = \big( \sqrt{a^2 + b^2} \big)^p \leq (|a| + |b|)^p \leq 2^{p-1}(|a|^p + |b|^p),\qquad p\geq 1.
	\end{align*}
\end{Proof}
\noindent By Proposition~\ref{tsconv},
\eqref{cpvr} defines a probability measure $\mu$
with $\mu(\ell^p(\N)) = 1$ 
 on the space $\ell^p(\N)$ of complex sequences.
 In Lemma~\ref{replem}, we determine the representation of the
 characteristic function of $\mu$ in the form \eqref{fjkldf34}.
\begin{lemma}\label{replem}
  Let $\alpha\in (0,1)$, $p\in [1,2]$, and let 
  $(k_n)_{n\geq 0} \in\ell^\alpha(\N)$ be a positive sequence. 
  Then, the distribution $\mu$ on $\ell^p(\N)$
                of the random
                series \eqref{cpvr} has characteristic functional \eqref{fjkldf34} with $R \equiv 0$, control measure
	        \begin{align}
                  \label{fkjhldf133} 
		\xi(dz) := \sum_{n=0}^\infty
		k_n^\alpha\big(
		\delta_{e_n} (dz) + \delta_{ie_n} (dz)
		\big), 
	\end{align}
	and L\'evy measures
	\begin{align*}
		\rho(e_n, du) = \rho(ie_n, du) =
		a_-\frac{e^{-\lambda_- |u| / k_n }}{|u|^{1+\alpha}} \mathbf{1}_{\R^-}(u) \du
		+
		a_+\frac{e^{-\lambda_+ u / k_n }}{u^{1+\alpha}} \mathbf{1}_{\R^+}(u) \du.
	\end{align*}
\end{lemma}
\begin{Proof}
  From \cite[Lemma 4.1]{TSDP},
  since $\theta$ is tempered stable on $\R$ with L\'evy measure \eqref{lefy},
  $k_n\theta$ is tempered stable on $\R$ with L\'evy measure
	\begin{align*}
	  \lambda_n(du) := k_n^\alpha 
           a_-\frac{e^{- \lambda_- |u| / k_n}}{|u|^{1+\alpha}} \mathbf{1}_{\R^-}(u) \du
	   +
           k_n^\alpha
			a_+\frac{e^{- \lambda_+ u / k_n }}{u^{1+\alpha}} \mathbf{1}_{\R^+}(u) \du.
	\end{align*}
	Next, 
        by independence of the sequences
        $(\theta_{1,n})_{n\geq 0}$,
        $(\theta_{2,n})_{n\geq 0}$,
 we have
	\begin{align*}
		&\int_{\ell^p(\N)} \exp \bigbracket{i \re\biganglebracket{z, x^*}} \mu(dz) = \E\bigsquarebracket{\exp \bigbracket{i \re\biganglebracket{\sum_{n=0}^\infty k_n(\theta_{1,n} + i\theta_{2,n})e_n, x^*}}}\\
	  &\quad = \prod_{n=0}^\infty
          \left(
          \E\bigsquarebracket{\exp \bigbracket{i \re\biganglebracket{ k_n\theta_{1,n}e_n, x^*}}}\E\bigsquarebracket{\exp \bigbracket{i \re\biganglebracket{ k_n\theta_{2,n}ie_n, x^*}}}
          \right)
          \\
	  &\quad = \prod_{n=0}^\infty
          \left(
          \E\bigsquarebracket{\exp \bigbracket{i \re\biganglebracket{e_n, x^*}k_n\theta_{1,n}}}\E\bigsquarebracket{\exp \bigbracket{i \re\biganglebracket{ie_n, x^*}k_n\theta_{2,n}}}
          \right)
          \\
	  &\quad = \prod_{n=0}^\infty
          \left(
          \exp\bigbracket{k_n^\alpha \int_{\R} \big(e^{iu\re\biganglebracket{e_n, x^*}}  - 1 - iu\kappa (u) \re\biganglebracket{e_n, x^*}\big) \frac{\lambda_n(du)}{k_n^\alpha}}
          \right. 
            \\
	    &\qquad\qquad
            \left.
            \times \exp\bigbracket{k_n^\alpha \int_{\R} \big(e^{iu\re\biganglebracket{ie_n, x^*}}  - 1 - iu\kappa (u) \re\biganglebracket{ie_n, x^*}\big) \frac{\lambda_n(du)}{k_n^\alpha}}
            \right)
              \\
	      &\quad = \exp\bigbracket{\int_{\ell^p(\N)}\int_{\R} \big(e^{iu\re\biganglebracket{z, x^*}}  - 1 - iu\kappa (u) \re\biganglebracket{z, x^*}\big)
                \sum_{n=0}^\infty
                \lambda_n(du) (\delta_{e_n}(dz) + \delta_{ie_n}(dz))}, 
	\end{align*}
	which is in the form \eqref{fjkldf34} with
        $
        \rho (e_n ,du) =\rho (ie_n ,du) =
        k_n^{-\alpha} \lambda_n(du)$, $n \geq 0$,
        and $\xi$ given by \eqref{fkjhldf133}, which
        is finite since $(k_n)_{n\geq 0} \in\ell^\alpha(\N)$.
\end{Proof}
\noindent
 We now present codifference bounds in the tempered stable setting.
\begin{lemma}
		\label{tempstabby}
	Let $\alpha\in (0,1)$ and suppose that $(k_n)_{n\geq 0} \in\ell^\alpha(\N)$ is a positive sequence. Let $p\in [1,2]$. Let $\mu$ be the distribution of \eqref{cpvr} on $\ell^p(\N)$. Then, the codifference bounds
	\begin{align}\label{whiffy1}
				\bigabs{C_\mu^= (x^*, y^*)} &\leq
		2^{4-p}(a_- \lambda_-^{\alpha-p} + a_+ \lambda_+^{\alpha-p})\Gamma\bigbracket{p-\alpha}\\
                \nonumber
		&\quad 
		\times \bigbracket{\sum_{n=0}^\infty
			k_n^{p}|\re\biganglebracket{e_n,x^*}\re\biganglebracket{e_n,y^*}|^{p/2} + \sum_{n=0}^\infty
			k_n^{p}|\imag\biganglebracket{e_n,x^*}\imag\biganglebracket{e_n,y^*}|^{p/2}}
	\end{align}
	and
	\begin{align}
          \label{whiffy2}
				\bigabs{C_\mu^{\neq} (x^*, y^*)} &\leq 
		2^{4-p}(a_- \lambda_-^{\alpha-p} + a_+ \lambda_+^{\alpha-p})\Gamma\bigbracket{p-\alpha}\\
                	\nonumber
		&\quad 
		\times \bigbracket{\sum_{n=0}^\infty
			k_n^{p}|\re\biganglebracket{e_n,x^*}\imag\biganglebracket{e_n,y^*}|^{p/2} + \sum_{n=0}^\infty
			k_n^{p}|\imag\biganglebracket{e_n,x^*}\re\biganglebracket{e_n,y^*}|^{p/2}}
	\end{align}
	hold for any $x^*, y^*\in (\ell^p(\N))^*$.
\end{lemma}

\begin{Proof}
	Using the control measure representation from Lemma~\ref{replem}, Lemma~\ref{fkjdf43} gives that for any $c > 0$ we have
	\begin{align}
          \label{whifffffff}
			\bigabs{C_\mu^= (x^*, y^*)} &\leq 2^{4-p} \underbrace{\sum_{n=0}^\infty k_n^\alpha |\re\biganglebracket{e_n,x^*}\re\biganglebracket{e_n,y^*}|^{p/2} I_c (e_n)}_{S_1(c)} +
		16\underbrace{\sum_{n=0}^\infty k_n^\alpha I_{c}' (e_n)}_{S_2(c)}
                \\
		\nonumber
                &\quad + 2^{4-p} \underbrace{\sum_{n=0}^\infty k_n^\alpha |\imag\biganglebracket{e_n,x^*}\imag\biganglebracket{e_n,y^*}|^{p/2} I_c (ie_n)}_{S_3(c)} +
		16\underbrace{\sum_{n=0}^\infty k_n^\alpha I_{c}' (ie_n)}_{S_2(c)},
	\end{align}
	where
	\begin{align*}
		I_c (e_n) = I_c(ie_n) & = 
		a_-\int_{-c}^0 \frac{ e^{- \lambda_- |u| / k_n}}{|u|^{1+\alpha-p}} \du
		+ a_+
		\int_0^c \frac{e^{- \lambda_+ u k_n}}{u^{1+\alpha-p}} \du
                \\
                  &\leq
		a_-\int_{-\infty}^0 \frac{ e^{-\lambda_- |u| / k_n}}{|u|^{1+\alpha-p}} \du
		+ a_+
		\int_0^\infty \frac{e^{-\lambda_+ u / k_n}}{u^{1+\alpha-p}} \du
		\\
		&=
		a_-\bigbracket{\frac{\lambda_-}{k_n}}^{\alpha-p}\int_{0}^\infty \frac{ e^{-x}}{x^{1+\alpha-p}} \dx
		+ a_+ \bigbracket{\frac{\lambda_+}{k_n}}^{\alpha-p}
		\int_0^\infty \frac{e^{-x}}{x^{1+\alpha-p}} \dx
		\\
		&= \frac{a_- \lambda_-^{\alpha-p} + a_+ \lambda_+^{\alpha-p} }{k_n^{\alpha-p}}\Gamma\bigbracket{p-\alpha},
\end{align*}
	and
        \begin{align*} 
		I_{c}' (e_n) = I'_c(ie_n) & =
		a_-\int_{-\infty}^{-c} \frac{ e^{-\lambda_- |u| /k_n}}{|u|^{1+\alpha}} \du
		+
		a_+\int_c^\infty \frac{ e^{-\lambda_+ u / k_n}}{u^{1+\alpha}} \du.
                \\
                &\leq 2(a_- + a_+)\int_c^\infty \frac{1}{u^{1+\alpha}} \du\\
                &= 2\frac{a_- + a_+}{\alpha c^\alpha}.
        \end{align*}
 This gives the bounds 
	\begin{align*}
	&	S_1(c) \leq (a_- \lambda_-^{\alpha-p} + a_+ \lambda_+^{\alpha-p} )\Gamma\bigbracket{p-\alpha}\sum_{n=0}^\infty k_n^{p}|\re\biganglebracket{e_n,x^*}\re\biganglebracket{e_n,y^*}|^{p/2}, 
	  \\
          & 
		S_3(c) \leq (a_- \lambda_-^{\alpha-p} + a_+ \lambda_+^{\alpha-p} )\Gamma\bigbracket{p-\alpha}\sum_{n=0}^\infty k_n^{p}|\imag\biganglebracket{e_n,x^*}\imag\biganglebracket{e_n,y^*}|^{p/2},
	\end{align*}
        independently of $c>0$, and 
	\begin{align*}
		S_2(c) 
	        \leq
                2\frac{a_- + a_+}{\alpha c^\alpha}
                \sum_{n=0}^\infty k_n^\alpha < \infty.
	\end{align*}
	Now taking the limit as $c\to\infty$ in \eqref{whifffffff}
         yields \eqref{whiffy1}. The proof for \eqref{whiffy2} is similar.
\end{Proof}
\noindent 
 We are now in a position to provide a class of tempered stable measures
 which admits a mixing operator, and to determine its mixing rate. 
\begin{prop}
	\label{exampprop}
	Let $\alpha\in (0,1)$, $p\in [1,2]$,
        and let $(\omega_n)_{n\geq 0}$ be a bounded positive weight sequence.
	In the framework of \eqref{cpvr},
	assume that the sequence $(k_n)_{n\geq 0}$ defined by $k_0 > 0$ and
	$$ 
	k_n :=
	\displaystyle
	k_0 \prod_{l=0}^{n-1}
	\frac{1}{\omega_l},
	\quad n \geq 1, 
	$$
	satisfies $(k_n)_{n\geq 0} \in \ell^\alpha(\N)$.
	Let $\mu$ be the distribution of \eqref{cpvr} on the Banach space $\ell^p(\N)$, and 
	consider the bounded weighted backward shift operator on $\ell^p(\N)$ defined by 
		\begin{align*}
		Te_0 := 0, \quad Te_{n+1} :=
		w_n 
		e_n,\qquad n\geq 0.
	\end{align*}
	The following are true.
	\begin{enumerate}
		\item $T$ admits $\mu$ as invariant measure.
		\item We have 
				\begin{align}
			\label{fjkdlsf2} 
			\sup_{x^*,y^*\in (\ell^p(\N))^*\setminus \{0\}}
			\frac{
				\bigabs{C_\mu^{=,\neq} (x^*, T^{*n}y^*)}}
			{                        \|x^*\|^{p/2} \|y^*\|^{p/2}}
			 \leq 
			2^{5-p}(a_- \lambda_-^{\alpha-p} + a_+ \lambda_+^{\alpha-p})\Gamma\bigbracket{p-\alpha}
			\sum_{l=0}^\infty k_l^{p / 2} k_{l+n}^{p / 2},
		\end{align}
			$n \geq 0$.
		In particular, $T$ is mixing
		by Theorem~\ref{fkds3},             provided that
		$\displaystyle \lim\limits_{n\to\infty}
		\sum_{l=0}^\infty k_l^{p / 2} k_{l+n}^{p / 2}  = 0$.
	\end{enumerate}
\end{prop}
\begin{Proof}
  \noindent
  1. 
  Let the random variable $X$ 
  be represented as in \eqref{cpvr}. Then,
		\begin{align*}
			TX &= \sum_{n=0}^\infty
			k_n\big(
			\theta_{1,n} + i\theta_{2,n}
			\big)Te_n\\
			&= 
			\sum_{n=1}^\infty
			k_n\omega_{n-1} \big(
			\theta_{1,n} + i\theta_{2,n}
			\big)e_{n-1}\\
			&= 
			\sum_{n=0}^\infty
			k_n \big(
			\theta_{1,n+1} + i\theta_{2,n+1}
			\big)e_n\\
			&\overset{d}{=} 
			\sum_{n=0}^\infty
			k_n \big(
			\theta_{1,n} + i\theta_{2,n}
			\big)e_n\\
			&= X,
		\end{align*}
		where the equality in distribution follows since the $\theta_{1,n}$'s and $\theta_{2,n}$'s are independent identical copies of the same random variable. The distribution of $TX$ is $\mu$ also, thus $T$ admits $\mu$ as invariant measure.
		
                \medskip

                \noindent
                2. We may represent $\mu$ with characteristic functional \eqref{fjkldf34} by Lemma~\ref{replem}, and thus the bounds of Lemma~\ref{tempstabby} apply. For $C_\mu^{=}(x^*, T^{*n}y^*)$, we have
		\begin{align*}
			\sum_{l=0}^\infty
			k_l^{p}|\re\biganglebracket{e_n,x^*}\re\biganglebracket{e_n,T^{*n}y^*}|^{p/2} & = 
			\sum_{l=0}^\infty k_{l+n}^{p} \bigabs{\re\biganglebracket{e_{l+n},x^*}\re\biganglebracket{e_l,y^*}\omega_l\dots\omega_{l+n-1}}^{p/2}\\
			&\quad \leq
			\|x^*\|^{p/2} \|y^*\|^{p/2} \sum_{l=0}^\infty k_{l+n}^{p}
			\prod_{j=l}^{l+n-1} \omega_j^{p/2}\\
			&\quad = \|x^*\|^{p/2} \|y^*\|^{p/2} \sum_{l=0}^\infty k_{l+n}^{p}
			\prod_{j=l}^{l+n-1} \bigbracket{\frac{k_j}{k_{j+1}}}^{p / 2}\\
			&\quad= \|x^*\|^{p/2} \|y^*\|^{p/2} \sum_{l=0}^\infty k_{l+n}^{p} \bigbracket{\frac{k_l}{k_{l+n}}}^{p / 2}\\
			&\quad= \|x^*\|^{p/2} \|y^*\|^{p/2}
			\sum_{l=0}^\infty k_l^{p / 2} k_{l+n}^{p / 2},
		\end{align*}
		and likewise
		\begin{align*}
			\sum_{l=0}^\infty
			k_l^{p}|\imag\biganglebracket{e_n,x^*}\imag\biganglebracket{e_n,T^{*n}y^*}|^{p/2} \leq \|x^*\|^{p/2} \|y^*\|^{p/2}
			\sum_{l=0}^\infty k_l^{p / 2} k_{l+n}^{p / 2}.
		\end{align*}
		Similarly, the same bounds hold for $C_\mu^{\neq}(x^*, T^{*n}y^*)$. We conclude by
		Lemma~\ref{tempstabby}. In particular,
			by Theorem~\ref{fkds3}, 
			$T$ is mixing on $\ell^p(\N)$ 
			when \eqref{fjkdlsf2} goes to zero
			as $n$ tends to infinity.
\end{Proof} 
	As with Proposition~\ref{extexp}, Lemma~\ref{replem}, and thus Proposition~\ref{exampprop}, can be applied to the symmetrized control measure
	\begin{align*}
		\xi(dz) := \frac{1}{2}\sum_{n=0}^\infty
		k_n^\alpha
		\big(
		\delta_{e_n} (dz) + \delta_{-e_n} (dz) + \delta_{ie_n} (dz) + \delta_{-ie_n} (dz)
		\big).
	\end{align*}
\noindent By a control of the quantity 
\begin{align*}
	\sum_{j=0}^\infty k_j^{p / 2} k_{j+l+n}^{p / 2}
	= 
	\sum_{j=0}^\infty 
	\frac{1}{
		(j+l+n+1)^{ \gamma /2}
		(j+1)^{\gamma /2}}
\end{align*}
 as in \eqref{b1}-\eqref{b2-2} of Example \ref{poicor},
 Proposition~\ref{exampprop} yields
 the following result. 
\begin{example}\label{exampcor}
	Let $\alpha \in (0, 1)$,
	$p\in [1,2]$, and $\gamma > 1$. In the context of Proposition~\ref{exampprop}, let
	$T$ be the bounded weighted backward shift operator on $E =\ell^p(\N)$ defined as 
	\begin{align*}
		Te_0 := 0, \quad Te_{n+1} :=
		\Bigg( 1 + \frac{1}{n+1}\Bigg)^{\gamma /p}
		e_n,\qquad n\geq 0,
	\end{align*}
	i.e. the weight sequence is
        $\omega_n = ( (n+2)/(n+1) )^{\gamma /p}$,
        $n\geq 0$. Define the coefficients of \eqref{cpvr} by
	\begin{align*}
		k_0 > 0,\qquad k_n = k_0\prod_{l=0}^{n-1}\frac{1}{\omega_l},
	\end{align*}
	and denote by $\mu$ the distribution of \eqref{cpvr}.
	Then $T$ admits $\mu$ as invariant measure,
	and $T$ is mixing with respect to $\mu$, with the rate 
\begin{flalign*}
    &\sup_{x^*,y^*\in \ell^p(\N))^*\setminus \{0\}}
    \frac{
        \bigabs{C_\mu^= (x^*, T^{*n}y^*)}}
    {\|x^*\|^{p/2}\|y^*\|^{p/2}}&
\end{flalign*}
\begin{empheq}[left={\displaystyle
      \leq \empheqlbrace}]{align}
\nonumber 
&
			\displaystyle 2^{5-p}(a_- \lambda_-^{\alpha-p} + a_+ \lambda_+^{\alpha-p})\Gamma\bigbracket{p-\alpha}k_0^p
			{\rm B} \left( 1 - \frac{\gamma}{2} , \gamma - 1 \right)
			n^{-(\gamma-1)}, \quad \quad \raisebox{0.0ex}{$\displaystyle
			1 < \gamma < 2,$}
\\[0.5ex]
\label{jkldsf11-1} 
&
\raisebox{0.5ex}{$\displaystyle 2^{5-p}(a_- \lambda_-^{\alpha-p} + a_+ \lambda_+^{\alpha-p})\Gamma\bigbracket{p-\alpha}k_0^p
			{\rm B} \left(\frac{\eps}{2} ,
			1 -\eps \right)n^{-(1-\eps)},$} \ \! \quad\quad\quad\quad \raisebox{0.5ex}{$\displaystyle\gamma \geq 2,$}
                        \end{empheq}
$n\geq 1$, for any $\eps \in (0,1)$ in \eqref{jkldsf11-1}. 
\end{example}
\noindent 
We now obtain convergence rates for more general functions in
the tempered stable setting.
\begin{corollary} 
	\label{fjkldf14}
	In the context of Proposition~\ref{exampprop},
	let $(a_j)_{j\in \N}$ and $(b_l)_{l\in \N}$
	be two complex $\ell^1(\N)$ sequences such that
	\begin{align*}
		\sum_{j=0}^\infty |a_{j}|\|T\|^{jp/2} < \infty.
	\end{align*}
	For any pair $(\phi,\psi)$ 
	chosen in Table~\ref{fctable},
	the functions 
	\begin{align}
		\label{genfg}
		f (z) := \sum_{j=0}^\infty a_j e^{i\phi
			( \langle z,T^{*j}x^*\rangle )}\qquad\text{and}\qquad g (z) := \sum_{l=0}^\infty b_l e^{i\psi
			( \langle z,T^{*l}y^* \rangle )}
	\end{align}
	are well-defined in $L^2(E,\mu)$, and if $T$ is mixing 
	we have 
	\begin{equation}
		\label{rate2} 
		\bigabs{I_n\bigbracket{f,g}}
		=
		O_{f,g,\mu ,T} \Bigg(
		\sum_{l=0}^\infty |b_l|
		\sum_{j=0 }^\infty k_j^{p / 2} k_{j+l+n}^{p / 2} 
		\Bigg), \qquad n \geq 0. 
	\end{equation} 
	\end{corollary}
\begin{Proof}
	As in Corollary~\ref{infseriesrate}, we have
	$f,g\in L^2(E,\mu)$. 
	We have the estimate
	\begin{align*}
		&\sum_{i=0}^\infty
		k_i^{p}|\re\biganglebracket{e_i,T^{*j}x^*}\re\biganglebracket{e_i,T^{*l}y^*}|^{p/2}\\
		&\quad = \sum_{i=\max ( 0, j-l ) }^\infty k_{i+l}^{p} \bigabs{\re\biganglebracket{e_{i+l-j},x^*}
			\omega_{i+l-j} 
			\dots
			\omega_{i+l-1}\re\biganglebracket{e_i, y^*}
			\omega_i
			\dots
			\omega_{i+l-1}}^{p/2}\\
		&\quad \leq
		\|T\|^{jp/2}
		\|x^*\|^{p/2}\|y^*\|^{p/2}
		\sum_{i=\max ( 0,j-l )}^\infty k_{i+l}^{p} \prod_{s=i}^{i+l-1}\omega_s^{p/2}
		\\
		& 
		\quad = \|T\|^{jp/2}
		\|x^*\|^{p/2}\|y^*\|^{p/2}
		\sum_{i=\max ( 0,j-l ) }^\infty k_i^{p / 2} k_{i+l}^{p/ 2} 
		,
		\quad
		j,l\geq 0. 
	\end{align*}
	The same estimate holds if we replace one or both of the real parts with the imaginary part.
		Hence, letting
	$$
	f_{j,x^*}(z) := e^{i\phi (\langle z,T^{*j}x^* \rangle )}, \quad 
	g_{l,y^*}(z)
	:= 
	e^{i\psi ( \langle z,T^{*l}y^* \rangle )},
	\quad
	j,l\geq 0,
	$$
	by Lemma~\ref{tempstabby}                   we have 
	\begin{align*} 
		& \bigabs{I_n\bigbracket{f_{j,x^*},g_{l,y^*}}}\\
		&\quad \leq  
		|C_\mu^{\phi, \psi }(T^{*j}x^*, T^{*(n+l)}y^*)|
		\exp\left(
		\max\left(
		\max_{n\geq 0} \big|C_\mu^{\phi , \psi } \big(T^{*n}x^*, y^*\big)\big|
		,
		\max_{n\geq 0} \big|C_\mu^{\phi , \psi } \big(x^*, T^{*n}y^*\big)\big|
		\right) \right) 
		\\
		&
		\quad
		=
		\|T\|^{jp/2}
				O_{x^*,y^*,\mu ,T} \Bigg(
		\sum_{i=\max(0,j-l-n)}^\infty k_i^{p / 2} k_{i+l+n}^{p/ 2} 
		\Bigg)
	\end{align*}
	and 
	\begin{align*}
		\nonumber 
		|I_n(f_{x^*},g_{y^*})| &\leq \sum_{j,l=0}^\infty |a_j||b_l| |I_n(f_j, g_l)|\nonumber\\
		& =  
				O_{x^*,y^*,\mu ,T}\Bigg(
		\sum_{j,l=0}^\infty 
		|a_j||b_l| \|T\|^{jp/2}
		\sum_{i=\max ( 0,j-l-n ) }^\infty
		k_i^{p / 2} k_{i+l+n}^{p / 2} 
		\Bigg),
	\end{align*}
	and \eqref{rate2} follows from the series convergence assumption.
\end{Proof}

\noindent Once again, by controlling the quantity 
\begin{align*}
	\sum_{j=0}^\infty k_j^{p / 2} k_{j+l+n}^{p / 2}
	= 
	\sum_{j=0}^\infty 
	\frac{1}{
		(j+l+n+1)^{ \gamma /2}
		(j+1)^{\gamma /2}}, 
\end{align*}
as in \eqref{b1}-\eqref{b2-2} of Example \ref{poicor},
 we obtain the following
from 
Corollary~\ref{fjkldf14}.
\begin{example}
	In the context of
	Example~\ref{exampcor}, if
	$$
	\sum_{j=0}^\infty |a_j|\|T\|^{jp/2} < \infty,
	$$ 	
	then
	$T$ admits the mixing rate
\begin{empheq}[left={\displaystyle
      \bigabs{I_n\bigbracket{f,g}} = 
      \empheqlbrace}]{align}
\nonumber 
       &
\displaystyle
			O_{f,g,\mu ,T} \big( n^{-(\gamma -1)}\big), \quad \displaystyle
			1 < \gamma < 2, 
\\[0.5ex]
\label{jkldsf11-2} 
&
\displaystyle
\raisebox{0.5ex}{$O_{f,g,\mu ,T} \big( n^{-(1-\eps ) }\big), \quad  \displaystyle  \gamma \geq 2$,}  
\end{empheq}
		$n\geq 1$,
	for any $\eps \in (0,1)$ in \eqref{jkldsf11-2}, 
		where $f,g$ are functions of the form \eqref{genfg}.
\end{example}
\noindent Finally, we also note that the rate
\eqref{rate2} holds for mixed exponential functions of the form 
\begin{align*}
	f(z) := \sum_{\phi\in\Phi}\sum_{j=0}^\infty a_{\phi , j} e^{i\phi(
		\langle z,T^{*j}x^* \rangle)}\qquad \text{and}\qquad 
	g(z) := \sum_{\psi\in\Psi}\sum_{l=0}^\infty b_{\psi , l}
	e^{i\psi (\langle z,T^{*l}y^* \rangle)}, 
\end{align*}
where $\Phi,\Psi$ 
are disjoint non-empty subsets
of $\{\pm\re(\cdot), \pm\imag(\cdot)\}$,  
and 
$(a_{\phi,j})_{j\in \N, \phi\in\Phi }$, 
$(b_{\psi,l})_{l\in \N , \psi\in\Psi }$
are complex $\ell^1(\N)$ sequences such that
\begin{align*}
	\sum_{j=0}^\infty |a_{\phi,j}|\|T\|^{jp/2} < \infty, \quad \phi \in \Phi.
\end{align*}
          
\footnotesize

\def\cprime{$'$} \def\polhk#1{\setbox0=\hbox{#1}{\ooalign{\hidewidth
  \lower1.5ex\hbox{`}\hidewidth\crcr\unhbox0}}}
  \def\polhk#1{\setbox0=\hbox{#1}{\ooalign{\hidewidth
  \lower1.5ex\hbox{`}\hidewidth\crcr\unhbox0}}} \def\cprime{$'$}


\begin{thebibliography}{Woy19}

\bibitem[App09]{applebk2}
D.~Applebaum.
\newblock {\em L\'{e}vy processes and stochastic calculus}, volume 116 of {\em
  Cambridge Studies in Advanced Mathematics}.
\newblock Cambridge University Press, Cambridge, second edition, 2009.

\bibitem[Bay15]{CLT}
F.~Bayart.
\newblock Central limit theorems in linear dynamics.
\newblock {\em Ann. Inst. Henri Poincar\'{e} Probab. Stat.}, 51(3):1131--1158,
  2015.

\bibitem[BM09]{DLO}
F.~Bayart and {\'{E}.}~Matheron.
\newblock {\em Dynamics of linear operators}, volume 179 of {\em Cambridge
  Tracts in Mathematics}.
\newblock Cambridge University Press, Cambridge, 2009.

\bibitem[CFS82]{cornfeld}
I.~P. Cornfeld, S.~V. Fomin, and Ya.~G. Sina{\u\i}.
\newblock {\em Ergodic theory}, volume 245 of {\em Grundlehren der
  Mathematischen Wissenschaften}.
\newblock Springer-Verlag, New York, 1982.

\bibitem[Dev13]{Devinck}
V.~Devinck.
\newblock Strongly mixing operators on {H}ilbert spaces and speed of mixing.
\newblock {\em Proceedings of the London Mathematical Society},
  106(3):1394--1434, 2013.

\bibitem[Dur10]{durrett2005}
R.~Durrett.
\newblock {\em Probability: Theory and Examples}.
\newblock Cambridge University Press, fourth edition, 2010.

\bibitem[FS13]{MCMIDP}
F.~Fuchs and R.~Stelzer.
\newblock Mixing conditions for multivariate infinitely divisible processes
  with an application to mixed moving averages and the sup{OU} stochastic
  volatility model.
\newblock {\em ESAIM: Probability and Statistics}, 17:455--471, 2013.

\bibitem[Kop95]{koponen}
I.~Koponen.
\newblock Analytic approach to the problem of convergence of truncated
  {L}\'{e}vy flights towards the {G}aussian stochastic process.
\newblock {\em Phys. Rev. E}, 52:1197--1199, 1995.

\bibitem[KT13]{TSDP}
U.~K\"uchler and S.~Tappe.
\newblock Tempered stable distributions and processes.
\newblock {\em Stochastic Processes and their Applications},
  123(12):4256--4293, 2013.

\bibitem[Lin86]{lindebook}
W.~Linde.
\newblock {\em Probability in {B}anach Spaces -- Stable and Infinitely
  Divisible Distributions}.
\newblock John Wiley \& Sons Ltd., 1986.

\bibitem[LT91]{LT}
M.~Ledoux and M.~Talagrand.
\newblock {\em Probability in {B}anach spaces}.
\newblock Springer-Verlag, 1991.

\bibitem[Mar70]{maruyama}
G.~Maruyama.
\newblock Infinitely divisible processes.
\newblock {\em Theory of Probability and Applications}, XV:1--22, 1970.

\bibitem[MP24]{MP23}
C.~Mau and N.~Privault.
\newblock Mixing of linear operators under infinitely divisible measures on
  {B}anach spaces.
\newblock {\em Journal of Mathematical Analysis and Applications}, 535:Paper No
  128160, 2024.

\bibitem[PP16]{Pitmanarticle}
E.J.G. Pitman and J.~Pitman.
\newblock A direct approach to the stable distributions.
\newblock {\em Advances in Applied Probability}, 48:261--282, 2016.

\bibitem[PV19]{MPMIDRF}
R.~Passeggeri and A.E.D. Veraart.
\newblock Mixing properties of multivariate infinitely divisible random fields.
\newblock {\em Journal of Theoretical Probability}, 32:1845--1879, 2019.

\bibitem[Ros87]{Rosinskialt}
J.~Rosi\'{n}ski.
\newblock Bilinear random integrals.
\newblock {\em Dissertationes Math. (Rozprawy Mat.)}, 259:71, 1987.

\bibitem[RZ96]{SCMIDP}
J.~Rosi\'nski and T.~Zak.
\newblock Simple conditions for mixing of infinitely divisible processes.
\newblock {\em Stochastic Processes and their Applications}, 61:277--288, 1996.

\bibitem[RZ97]{EEWMIDP}
J.~Rosi\'nski and T.~Zak.
\newblock Equivalence of ergodicity and weak mixing for infinitely divisible
  processes.
\newblock {\em Journal of Theoretical Probability}, 10(1):73--86, 1997.

\bibitem[Sat99]{sato}
K.~Sato.
\newblock {\em L\'evy processes and infinitely divisible distributions},
  volume~68 of {\em Cambridge Studies in Advanced Mathematics}.
\newblock Cambridge University Press, Cambridge, 1999.

\bibitem[Sch70]{Owls}
L.~Schwartz.
\newblock Les applications {$O$}-radonifiantes dans les espaces de suites.
\newblock {\em S\'eminaire d'Analyse fonctionnelle (dit ``Maurey-Schwartz''),
  exp. no 26}, pages 1--19, 1969--1970.
\newblock Ecole Polytechnique, Centre de Math\'ematiques.

\bibitem[Tor77]{Tortrat}
A.~Tortrat.
\newblock Lois {e$(\lambda )$} dans les espaces vectoriels et lois stables.
\newblock {\em Z. Wahrscheinlichkeitstheorie und Verw. Gebiete},
  37(2):175--182, 1976/77.

\bibitem[WA57]{akutowicz}
N.~Wiener and E.J. Akutowicz.
\newblock The definition and ergodic properties of the stochastic adjoint of a
  unitary transformation.
\newblock {\em Rend. Circ. Mat. Palermo (2)}, 6:205--217; addendum, 349, 1957.

\bibitem[Woy19]{GMBS}
W.A. Woyczy\'nski.
\newblock {\em Geometry and martingales in {B}anach spaces}.
\newblock CRC Press, Boca Raton, FL, 2019.

\end{thebibliography}
\end{document}